\def\yes{\if00}
\def\no{\if01}
\def\iftwelvept{\no}

\def\ifusepdf{\no}
\def\ifpsfont{\yes}

\iftwelvept
\documentclass[leqno,12pt]{amsart}
\else
\documentclass[leqno]{amsart}
\fi
\usepackage{amssymb}
\usepackage{amscd}
\usepackage{latexsym}
\usepackage{verbatim}
\usepackage[all]{xy}

\ifusepdf
\usepackage{hyperref}
\else\fi
\ifpsfont
\usepackage[T1]{fontenc}
\usepackage{times}
\else\fi

\iftwelvept
\else\fi

\newtheorem{Theorem}{Theorem}[section]
\newtheorem{Proposition}[Theorem]{Proposition}
\newtheorem{Lemma}[Theorem]{Lemma}
\newtheorem{Corollary}[Theorem]{Corollary}
\newtheorem{Claim}{Claim}[Theorem]

\theoremstyle{definition}
\newtheorem{Definition}[Theorem]{Definition}
\newtheorem{Remark}[Theorem]{Remark}

\renewcommand{\theTheorem}{\arabic{section}.\arabic{Theorem}}

\newcommand{\ZZ}{{\mathbb{Z}}}

\newcommand{\RR}{{\mathbb{R}}}
\newcommand{\CC}{{\mathbb{C}}}
\newcommand{\PP}{{\mathbb{P}}}

\newcommand{\Aff}{{\mathbb{A}}}

\newcommand{\Gal}{\operatorname{Gal}}

\newcommand{\Kcal}{{\mathcal{K}}}
\newcommand{{\Vcal}}{{\mathcal{V}}}

\newcommand{\Proof}{{\sl Proof.}\quad}

\newcommand{\QED}{{\unskip\nobreak\hfil\penalty50\quad\null\nobreak\hfil
{$\Box$}\parfillskip0pt\finalhyphendemerits0\par\medskip}}

\begin{document}

\title[Canonical heights for regular automorphisms]{Local and global canonical height functions for affine space regular automorphisms} 
\author{Shu Kawaguchi}
\address{Department of Mathematics, Graduate School of Science,
Osaka University, Toyonaka, Osaka 560-0043, Japan}
\email{kawaguch@math.sci.osaka-u.ac.jp}
\thanks{Supported in part by MEXT grant-in-aid for young scientists (B) 21740018}

\begin{abstract}
Let $f: \Aff^N \to \Aff^N$ be a regular polynomial automorphism
defined over a number field $K$. For each place $v$ of $K$, we
construct the $v$-adic Green functions $G_{f,v}$ and $G_{f^{-1},v}$ 
(i.e., the $v$-adic canonical height functions) for $f$ and $f^{-1}$. 
Next we introduce for $f$ the notion of good reduction at $v$, and 
using this notion, we show that the sum of $v$-adic Green functions 
over all $v$ gives rise to a canonical height function for $f$ that 
satisfies the Northcott-type finiteness property. 
Using \cite{Ka}, we recover results on arithmetic properties of 
$f$-periodic points and non $f$-periodic points. We also obtain an estimate 
of growth of heights under $f$ and $f^{-1}$, which is independently 
obtained by Lee by a different method. 
\end{abstract} 

\maketitle

\section*{Introduction}
\renewcommand{\theTheorem}{\Alph{Theorem}} 
Height functions are one of the basic tools in Diophantine geometry. On Abelian 
varieties defined over a number field, there exist N\'eron--Tate's canonical height functions 
that behave well relative to the $[n]$-th power map. 
Tate's elegant construction is via a global method 
using a relation of an ample divisor relative to the $[n]$-th power map. 
N\'eron's construction 
is via a local method, and gives deeper properties 
of the canonical height functions. 
Both constructions 
are useful in studying arithmetic properties of Abelian varieties. 

In \cite{Ka}, we showed the existence of canonical height functions for affine plane 
polynomial automorphisms of dynamical degree $\geq 2$. Our construction was 
via a global method using the effectiveness of a certain divisor on a certain rational surface. 
In this paper, we use a local method to 
construct a canonical height function for affine space regular automorphism 
$f: \Aff^N \to \Aff^N$, which coincides with the one in \cite{Ka} when $N=2$. 
We note that arithmetic properties of polynomial automorphisms over number fields
have been studied, for example, by Silverman \cite{SiHenon}, Denis \cite{De}, Marcello 
\cite{Ma1, Ma2} and the author \cite{Ka}. 

We recall the definition of regular polynomial automorphisms. 
Let $f: \Aff^N \to \Aff^N$ be a  polynomial automorphism of degree $d \geq 2$
defined over a field, and $\overline{f}: \PP^N\dasharrow\PP^N$ denote its birational 
extension to $\PP^N$. We write $f^{-1}$ for the inverse of $f$, $d_-$ for the degree 
of $f^{-1}$, and $\overline{f^{-1}}$ for its birational extension to $\PP^N$. 
Then $f$ is said to be {\em regular} if the intersection of the set of indeterminacy 
of $\overline{f}$ and that of  $\overline{f^{-1}}$ is empty over an algebraic closure of 
the field (cf. Definition~\ref{def:regular:autom} and Remark~\ref{rmk:regular:autom:def}). 
Over $\CC$, dynamical properties of affine space regular polynomial automorphisms $f$ 
are deeply studied, in which the Green function for $f$ plays a pivotal role (see \cite[\S2]{Sib}). 

In \S\ref{sec:local:theory} and \S\ref{sec:local:theory:II}, 
we construct a Green function (a local canonical height function) 
for $f$ over an algebraically closed field $\Omega$
with non-trivial non-archimedean absolute value $|\cdot|$. 
For $x = (x_1, \ldots, x_N) \in \Omega^N$, we set 
$\Vert x \Vert = \max_{1 \leq i \leq N}\{|x_i|\}$. 
Our results are put together as follows. 

\begin{Theorem}[cf. Proposition~\ref{prop:local:Green:function}, Lemma~\ref{lemma:G:above}, Theorem~\ref{thm:local:II}]
\label{thm:A}
Let $f: \Aff^N \to \Aff^N$ be a regular polynomial automorphism of 
degree $d \geq 2$ defined over $\Omega$. With the notation as above, 
we have the following. 
\begin{enumerate}
\item
For all $x \in \Aff^N(\Omega)$, the limits 
\[
  \lim_{n\to+\infty}\frac{1}{d^n} \log\max\{\Vert f^n(x) \Vert, 1\}
  \quad\text{and}\quad
  \lim_{n\to+\infty}\frac{1}{d_-^n}  \log\max\{\Vert f^{-n}(x) \Vert, 1\}
\]
exist. We respectively write $G_f(x)$ and $G_{f^{-1}}(x)$ 
for the limits, which we call Green functions for $f$ and $f^{-1}$.
\item
There are constants $c_f, c_{f^{-1}} \in \RR$ such that 
\[
  G_f(\cdot) \leq \log\max\{\Vert\cdot\Vert, 1\} + c_f
  \quad\text{and}\quad
  G_{f^{-1}}(\cdot) \leq \log\max\{\Vert\cdot\Vert, 1\} + c_{f^{-1}}
  \quad\text{on $\Aff^N(\Omega)$}.
\]  
\item
There are subsets $V^+, V^-$ of $\Aff^N(\Omega)$ with 
$V^+ \cup V^- = \Aff^N(\Omega)$  
and constants 
$c^+, c^- \in \RR$ such that 
\begin{align*}
  G_f(\cdot) & \geq \log\max\{\Vert\cdot\Vert, 1\} + c^+ \quad\text{on $V^+$}, \\
  G_{f^-}(\cdot) & \geq \log\max\{\Vert\cdot\Vert, 1\} + c^- \quad\text{on $V^-$}. \\
\end{align*}
\end{enumerate}
\end{Theorem}

Over $\CC$, Green functions are constructed using compactness 
arguments (cf. \cite[\S2]{Sib}). Here we use more algebraic arguments based on  
Hilbert's Nullstellensatz. Our construction of $V^\pm, c^{\pm}$ is rather delicate 
with a choice of two parameters $\varepsilon$ and $\delta$, which behaves well when we work 
over number fields in \S\ref{sec:global:theory} and \S\ref{sec:applications}. 
We note that, over $\CC$, our construction gives a different proof of the existence of 
Green functions with more explicit estimates (see \S\ref{sec:Green:over:C}). 
In \S\ref{sec:local:theory:III}, we continue to study some basic properties 
of regular polynomial automorphisms $f$ over $\Omega$, characterizing 
the set of escaping points by $G_f$ and showing a filtration 
property for $f$. 

\smallskip
Now we turn our attention to number fields. 
Let $f: \Aff^N \to \Aff^N$ be a  polynomial automorphism 
defined over a number field $K$. For each place $v$ of $K$, let $K_v$ denote 
the completion of $K$ with respect to $v$, and $\overline{K_v}$ an algebraic closure 
of $K_v$. Then $f$ induces a regular polynomial automorphism over 
$\overline{K_v}$, so we have Green functions $G_{f, v}, G_{f^{-1}, v}$ and 
estimates with $c_{f, v}, c_{f^{-1}, v}, c^{\pm}_{v}$ as in Theorem~\ref{thm:A}. (Here we use 
the suffix $v$ to indicate that we work over $\overline{K_v}$. See \S\ref{sec:Green:over:C} 
when $v$ is Archimedean.)   

We want to define the canonical height functions $\widehat{h}_f^+, \widehat{h}_f^-$ for $f$ 
as the sum of $G_{f, v}, G_{f^{-1}, v}$ over all the places $v$ of $K$. To this end, 
we introduce the notion of good 
reduction at an non-Archimedean place $v$ of $K$. 
Let $\overline{R_v}$ denote the ring of integers of $\overline{K_v}$, and $\widetilde{k_v}$ 
the residue field. Recall that the notion of good reduction for an endomorphism $\varphi$ of 
$\PP^1$ over $\overline{K_v}$ is introduced in Morton--Silverman \cite{MS}, which 
means that $\varphi$ extends to a morphism over $\overline{R_v}$ and the induced morphism 
$\widetilde{\varphi}$ over $\widetilde{k_v}$ has the same degree as $\varphi$. Here we say 
that a regular polynomial automorphism $f: \Aff^N\to\Aff^N$ has {\em good reduction} at $v$ 
if $f$ extends to an automorphism over $\overline{R_v}$ and the the induced morphism 
$\widetilde{f}$ over $\widetilde{k_v}$ is again a regular polynomial automorphism such 
that the degrees of $\widetilde{f}$ and $\widetilde{f}^{-1}$ are the same as the degrees of 
$f$ and $f^{-1}$ respectively (see Definition~\ref{def:good:reduction} for the precise definition). 

Using the notion of good reduction, we show the existence of canonical height functions. 
Let $h: \Aff^N(\overline{K}) \to \RR$ denote the usual logarithmic Weil height function. 

\begin{Theorem}[cf. Proposition~\ref{prop:good:reduction:are:ae}, Theorem~\ref{thm:main}] 
\label{thm:B}
Let $f: \Aff^N\to\Aff^N$ be a regular polynomial automorphism 
of degree $d \geq 2$ over a number field $K$. Let $d_- \geq 2$ 
denote the degree of $f^{-1}$. 
\begin{enumerate}
\item
$f$ has good reduction at $v$ except for finitely many places. Further, 
if this is the case, we can take the constants $c_{f,v} = c_{f^{-1}, v} = c^\pm_v =0$ 
in Theorem~\textup{\ref{thm:A}}. 
\item
For all $x \in \Aff^N(\overline{K})$, 
the limits 
\begin{equation}
\label{eqn:thm:B}
  \widehat{h}^+_f(x):= \lim_{n\to+\infty}\frac{1}{d^n} h(f^n(x)) 
  \quad\text{and}\quad
  \widehat{h}^-_{f}(x):= \lim_{n\to+\infty}\frac{1}{d_-^n} h(f^{-n}(x))
\end{equation}
exist. Further, we have the decomposition into the sum of 
local Green functions\textup{:} 
\[
  \widehat{h}^+_f(x) = \sum_{v \in M_K} n_v G_{f,v}(x) 
  \quad\text{and}\quad
  \widehat{h}^-_{f}(x) = \sum_{v \in M_K} n_v G_{f^{-1},v}(x). 
\]
\item
We define $\widehat{h}_f: \Aff^N(\overline{K})\to\RR$ by 
$\widehat{h}_f := \widehat{h}^+_f + \widehat{h}^-_{f}$. 
Then $\widehat{h}_f$ satisfies $\widehat{h}_f \gg\ll h$ 
and 
\[
\frac{1}{d} \widehat{h}_f\circ f + 
\frac{1}{d_-} \widehat{h}_f\circ f^{-1} 
= \left(1+\frac{1}{d d_-}\right) \widehat{h}_f.  
\] 
Further, 
for $x \in \Aff^N(\overline{K})$, we have 
\[
  \widehat{h}_f(x) = 0 
  \;\Longleftrightarrow\;
  \widehat{h}^+_f(x) = 0 
  \;\Longleftrightarrow\;
  \widehat{h}^-_f(x) = 0 
  \;\Longleftrightarrow\;
  \text{$x$ is $f$-periodic.}
\]
\end{enumerate}
\end{Theorem}

We note that in \cite{Ka}, we have defined $\widehat{h}^{+}_f(x)$ as  $\limsup_{n\to\infty}    
\frac{1}{d^n} h(f^n(x))$, and similarly for $\widehat{h}^{-}_f$. Theorem~\ref{thm:B} 
shows that $\{   
\frac{1}{d^n} h(f^n(x))\}_{n=0}^{+\infty}$ and $\{   
\frac{1}{d_-^n} h(f^{-n}(x))\}_{n=0}^{+\infty}$
are in fact convergent sequences, i.e., $\limsup$ can 
be replaced by $\lim$ in \eqref{eqn:thm:B}. 

Using estimates on local Green functions over all places, 
we obtain the following estimate on global canonical height 
functions for all $N \geq 2$ (see \cite[\S4]{Ka}, \cite[Conjecture~3]{Sil}, \cite[Conjecture~7.18]{SilBook}). 
This result is obtained by ChongGyu Lee \cite{Lee} independently. 
His proof uses a global method and is based 
on the effectiveness of a certain divisor (as done for $N=2$ in \cite{Ka}). 

\begin{Corollary}[cf. Theorem~\ref{thm:former:assumption}]
\label{cor:C}
Let $f: \Aff^N \to \Aff^N$ be a regular polynomial automorphism 
over a number field $K$. With the notation as above, there 
exists a constant $c \geq 0$ 
such that 
\begin{equation}
\label{eqn:cor:C} 
  \frac{1}{d}h(f(x)) + \frac{1}{d_-}h(f^{-1}(x))
  \geq \left(1 + \frac{1}{d d_-}\right) h(x) -c  
\end{equation}
for all $x \in \Aff^N(\overline{K})$.  Further, we have 
\[
  \liminf_{\substack{
  x \in \Aff^N(\overline{K})\\ 
  h(x) \to \infty}} 
  \frac{
  \frac{1}{d}h(f(x)) + \frac{1}{d_-}h(f^{-1}(x))}{h(x)} 
  = 1 + \frac{1}{d d_-}. 
\]
\end{Corollary}

Since \eqref{eqn:cor:C} holds, by the argument of \cite{Ka}, 
we recover the results on $f$-periodic points and 
refine the results on non $f$-periodic points in \cite{SiHenon, De, Ma1, Ma2}. 
For $x \in \Aff^N(\overline{K})$, let  
$O_f(x) := \{f^n(x) \mid n \in \ZZ\}$ denote the $f$-orbit of $x$. 
If $O_f(x)$ is infinite, we have the canonical 
height $\widehat{h}(O_f(x))$ of $O_f(x)$ (see Eqn.~\eqref{eqn:can:ht:orbit}). 

\begin{Corollary}[cf. Corollary~\ref{eqn:can:ht:orbit}, Corollary~\ref{cor:refined}]
\label{cor:D}
Let $f: \Aff^N \to \Aff^N$ be a regular polynomial automorphism 
over a number field $K$. With the notation as above, 
we have the following. 
\begin{enumerate}
\item
The set of $f$-periodic points in $\Aff^N(\overline{K})$ 
is a set of bounded height. 
\item
For any infinite orbit $O_f(x)$, 
\[
  \#\{y \in O_f(x) \mid h(y) \leq T\}
  = 
  \left(\frac{1}{\log d} + \frac{1}{\log d_-} \right)\log T 
  - \widehat{h}(O_f(x)) + O(1)
\]
as $T \to +\infty$. 
\end{enumerate}
\end{Corollary}

\smallskip
{\sl Acknowledgment.}\quad My deep thanks go to Professor
Joseph H. Silverman for his encouragement and discussions. 
Some preliminary work was done 
while I was staying at the IMJ in Paris and the CRM in Barcelona 
from October of 2005 to July of 2006. 
I thank Professors Tien-Cuong Dinh, Charles Favre 
and Nessim Sibony for helpful conversations, 
with whom I had opportunities to talk during my stay 
in Paris. 

\renewcommand{\theTheorem}{\arabic{section}.\arabic{Theorem}}
\setcounter{equation}{0}
\medskip
\section{Non-archimedean Green functions for polynomial maps} 
\label{sec:local:theory}
Let $\Omega$ be an algebraically closed field 
with non-trivial non-archimedean absolute value $|\cdot|$, and  
$R$ its ring of integers. For a point 
$x= (x_1, \ldots, x_N) \in \Aff^N(\Omega)$, the norm of $x$ is 
defined by $\Vert x \Vert = \max_{i = 1, \ldots, N} \{|x_i|\}$. 
As usual, we set 
$\log^+(a) := \log\max\{a, 1\}$ for $a \in \RR_{\geq 0}$, 
so that $\log^+\Vert x \Vert = \log\max\{\Vert x \Vert, 1\} 
= \log\Vert (x,1) \Vert$. 

Let $f = (f_1, \ldots, f_N): \Aff^N\to\Aff^N$ be 
a polynomial map of degree $d \geq 2$ defined over $\Omega$, 
where $f_1(X), \ldots, f_N(X)$ are polynomials in 
$\Omega[X_1, \ldots, X_N]$ such that  
$d = \max_{i = 1, \ldots, N}\{\deg f_i\}$. 
We write $F_i(X, T) := T^d f_i(X/T)\in \Omega[X_1, \ldots, X_N, T]$ for 
homogenization of $f_i$. 
Let $\overline{f} = (F_1: \cdots: F_N: T^d): 
\PP^N \dasharrow \PP^N$ denote the extension of $f$ 
to $\PP^N$. 
We put $F:=(F_1, \ldots, F_N, T^d): \Aff^{N+1} \to \Aff^{N+1}$, 
which is a lift of $\overline{f}$. 

For the composition $f^n = f\circ\cdots\circ f$, we write 
$f^n = (f_1^n, \ldots, f_N^n)$. Similarly, for the composition 
$F^n = F\circ\cdots\circ F$, we write $F^n = (F_1^n, \ldots, F_N^n, T^{d^n})$.
Let $d_n$ denote the degree of $f^n$, and let 
$F_{ni}(X,T)= T^{d_n} f_i^n(X/T)\in \Omega[X_1, \ldots, X_N, T]$ be 
homogenization of $f_i^n$. 
Since $F_i^n(X, 1) = f_i^n(X) = F_{ni}(X, 1)$, the degree counting gives 
$F_i^n(X,T) = T^{d^n - d_n} F_{ni}(X,T)$. 

\begin{Proposition}
\label{prop:local:Green:function}
Let $f: \Aff^N\to\Aff^N$ be 
a polynomial map of degree $d \geq 2$ defined over $\Omega$. 
Then, for all $x \in \Aff^N(\Omega)$, 
$\displaystyle{\frac{1}{d^n} \log^+\Vert f^n(x)\Vert}$ 
converges to a non-negative real number as $n \to +\infty$. 
\end{Proposition}

\Proof
We take an $r \in R$ so that $r F_i \in R[X, T]$ for 
all $i =1, \ldots, N$. 
We set 
\[
  a_n := \frac{1}{d^n} \log^+\Vert f^n(x)\Vert, \quad
  b_n := \frac{1}{d^n} \log\Vert F^n(x,1)\Vert, \quad
  c_n := \frac{1}{d^n} \log\Vert (rF)^n(x,1)\Vert, 
\]
where $rF = (r F_1, \ldots, rF_N, r T^d)$. 
We claim that 
{\allowdisplaybreaks
\begin{equation}
\label{eqn:local:Green:function}
  a_n = b_n = c_n - \frac{1-d^{-n}}{d-1} \log |r|.
\end{equation}
Indeed, the first equality follows from 
\begin{align*}
  a_n 
  & = \frac{1}{d^n} \log^+\Vert f^n(x) \Vert \\
  & = \frac{1}{d^n} \log \max\{
    |F_{n1}(x, 1)|, \ldots, |F_{nN}(x, 1)|, 1 \} \\
  & = \frac{1}{d^n} \log \max\{
    |F_{1}^n(x, 1)|, \ldots, |F_{N}^n(x, 1)|, 1 \} 
  = b_n.
\end{align*}
}
The second equality follows from 
$(rF)^n = r^{1+d +\cdots + d^{n-1}} F^n = r^{\frac{d^n -1}{d-1}} F^n$.  
It follows from $\Vert (rF)(x, 1)\Vert \leq \Vert (x, 1)\Vert^d$ that 
\[
  \frac{1}{d^n} \log\Vert (rF)^n(x,1)\Vert 
  \leq  \frac{1}{d^n} \log\Vert (rF)^{n-1}(x,1)\Vert^d 
  \leq \frac{1}{d^{n-1}} \log\Vert (rF)^{n-1}(x,1)\Vert.   
\]
In other words, $\{c_n\}_{n=1}^{+\infty}$ is a non-increasing sequence. 
Eqn.~\eqref{eqn:local:Green:function} implies that 
$\{c_n\}_{n=1}^{+\infty}$ is bounded from below. Indeed, 
since $a_n$ is nonnegative and $|r| \leq 1$, we have 
$c_n \geq a_n + \frac{1}{d-1} \log|r| \geq \frac{1}{d-1} \log|r|$. 
Thus $\lim_{n\to+\infty} c_n$ exists. 
Eqn.~\eqref{eqn:local:Green:function} then gives 
the existence of $\lim_{n\to+\infty} a_n$, 
which is nonnegative from the definition. 
\QED

Proposition~\ref{prop:local:Green:function} allows the following 
definition.  

\begin{Definition}
\label{def:Green:function}
For a polynomial map $f:\Aff^N\to\Aff^N$ defined over $\Omega$, 
we define the non-negative function 
$G_f: \Aff^N(\Omega) \to \RR$ by 
\[
  \displaystyle{G_f(x) := \lim_{n\to+\infty}\frac{1}{d^n} 
  \log^+\Vert f^n(x)\Vert}
  \qquad 
  \text{for $x \in \Aff^N(\Omega)$}, 
\] 
and call it the {\em Green function} for $f$.  
\end{Definition}

\begin{Lemma}
\label{lemma:G:above}
Let $C_f^\prime$ be the maximum of the absolute value of 
all the coefficients of $f_i(X) $ for $1 \leq i \leq N$, 
and we set $c_f = \frac{1}{d-1}\log \max\{C_f^\prime, 1\}$. Then 
\[
  G_f(\cdot) \leq \log^+\Vert \cdot \Vert + c_f
  \quad\text{on $\Aff^N(\Omega)$}. 
\]
\end{Lemma}

\Proof
We take $r \in R$ such that $|r| = \frac{1}{\max\{C_f^\prime, 1\}}$. Then 
$rF_i \in R[X,T]$ for all $i = 1, \ldots, N$. 
From the proof of Proposition~\ref{prop:local:Green:function}, 
we have
\[
  G_f(x) \leq \lim_{n\to+\infty} c_n - \frac{1}{d-1}\log|r| 
  \leq c_0 - \frac{1}{d-1}\log|r| 
  = \log^+\Vert x \Vert - \frac{1}{d-1}\log|r|. 
\]
Hence we get the assertion. 
\QED

Lemma~\ref{lemma:AS} below shows that, 
for some polynomial maps $f$, $G_f$ is not interesting. 
However, we will see in the next section that 
$G_f$ enjoys nice properties for regular polynomial 
automorphisms $f$ (see Definition~\ref{def:regular:autom} 
and Theorem~\ref{thm:local:II}). 

To state Lemma~\ref{lemma:AS}, 
we recall that a polynomial map $f$ is said to be {\em algebraically stable} 
if $d_n = d^n$ for all $n \geq 1$ (\cite[\S1.4]{Sib}). 

\begin{Lemma}
\label{lemma:AS}
If $f$ is not algebraically stable,  
then $G_f(x) = 0$ for all $x \in \Aff^N(\Omega)$. 
\end{Lemma}

\Proof
We take $n_0$ such that $d_{n_0} < d^{n_0}$, and we put 
$g = f^{n_0}$. Proposition~\ref{prop:local:Green:function} 
tells us that 
$\frac{1}{d_{n_0}^m} \log^+\Vert g^m(x)\Vert$
converges to a non-negative number as $m \to+\infty$. 
Hence
\[
  \frac{1}{d^{n_0 m}} \log^+\Vert f^{n_0 m}(x)\Vert
  = \left(\frac{d_{n_0}}{d^{n_0}}\right)^m 
  \frac{1}{d_{n_0}^m}\log^+\Vert g^m(x)\Vert
  \longrightarrow 0 
  \qquad \text{($m \to +\infty$)}.   
\]
From Proposition~\ref{prop:local:Green:function}, we get 
$G_f(x) = 0$. 
\QED

\medskip
\section{Non-archimedean Green functions for regular automorphisms}
\label{sec:local:theory:II}
\setcounter{equation}{0}
In this section we consider polynomial automorphisms. 
Let $f: \Aff^N\to\Aff^N$ be a polynomial automorphism
of degree $d \geq 2$ defined over an algebraically closed field $\Omega$ with 
non-trivial non-archimedean absolute value. 

As before, 
let $\overline{f} = (F_1(X,T):\cdots:F_N(X,T):T^d):
\PP^N\dasharrow\PP^N$ denote the extension of $f$ 
to $\PP^N$. 
For the inverse $f^{-1}: \Aff^N\to\Aff^N$ of $f$, 
we denote by $d_-$ the degree of $f^{-1}$. 
Note that the integer $d_- \geq 2$ may be 
different from $d$. We write 
$\overline{f^{-1}}=(G_1(X,T):\cdots:G_N(X,T):T^{d_-}): \PP^N\dasharrow\PP^N$ 
for the extension of $f^{-1}$ to $\PP^N$.  

Let $I_+$ and $I_-$ denote the set of indeterminacy of 
$\overline{f}$ and $\overline{f^{-1}}$, respectively: 
\begin{gather*}
I_+ = \{(x:0)\in \PP^N(\Omega) \mid F_1(x,0)=\cdots=F_N(x,0)=0\}, \\
I_- = \{(x:0)\in \PP^N(\Omega) \mid G_1(x,0)=\cdots=G_N(x,0)=0\}.
\end{gather*}

\begin{Definition}[see {\cite[\S2.2]{Sib}}]
\label{def:regular:autom}
A polynomial automorphism $f: \Aff^N\to\Aff^N$ is called {\em regular} 
if $I_+ \cap I_- = \emptyset$. 
\end{Definition}

\begin{Remark}
\label{rmk:regular:autom:def}
The definition of regular polynomial automorphisms works over  
any algebraically closed field. 
\end{Remark}

The purpose of this section is to prove the following theorem, 
which says that the Green functions for regular automorphisms 
exhibit nice properties. 

\begin{Theorem}
\label{thm:local:II}
Let $\Omega$ be an algebraically closed field 
with non-trivial non-archimedean valuation, and 
$f: \Aff^N\to\Aff^N$ a regular polynomial automorphism 
over $\Omega$. Then there are open subsets ${V}^+, {V}^-$ of 
$\Aff^N(\Omega)$ and constants $c^+, c^- \in \RR$ with the following 
properties. 
\begin{enumerate}
\item[(i)] $G_f(\cdot) \geq \log^+\Vert \cdot\Vert + c^+$ on ${V}^+$. 
\item[(ii)] $G_{f^{-1}}(\cdot) \geq \log^+\Vert \cdot\Vert + c^-$ on ${V}^-$. 
\item[(iii)] ${V}^+ \cup {V}^- = \Aff^N(\Omega)$. 
\end{enumerate}
\end{Theorem}

\begin{Remark}
Over $\CC$, corresponding results (and much more) were 
established by Sibony \cite[\S2.2]{Sib}. 
Here, since $\Aff^N(\Omega)$ is not locally compact in general, 
we give a different proof that is more algebraic in nature 
based on Hilbert's Nullstellensatz. 
We also give $V^+, V^-$ and $c^+, c^-$ with precise estimates, 
so that they work well when we introduce the notion of 
good reduction in \S\ref{sec:good:reduction}. 
\end{Remark}

Before proving Theorem~\ref{thm:local:II}, we show several 
lemmas.  We begin by introducing some notation. 
Since $I_+ \cap I_-$ is empty, 
$F_1(X,0), \ldots, F_N(X,0), G_1(X,0), \ldots, 
G_N(X,0)$ have no solutions in common other than $0$. 
Thus, for each $1 \leq i \leq N$, 
there are polynomials $P_{ij}(X), Q_{ij}(X) \in \Omega[X]$ 
for $1 \leq j \leq N$ such that 
\begin{equation}
\label{eqn:Pij:Qij}
  \sum_{j=1}^N P_{ij}(X)F_j(X,0) 
  + \sum_{j=1}^N Q_{ij}(X)G_j(X,0) 
  = X_i^m
\end{equation}
with some $m \geq 1$. Hence there is 
a polynomial $R_i(X, T)\in\Omega[X,T]$ such that 
\begin{equation}
\label{eqn:pqr}
  \sum_{j=1}^N P_{ij}(X)F_j(X,T) 
  + \sum_{j=1}^N Q_{ij}(X)G_j(X,T) + T R_i(X, T)
  = X_i^m. 
\end{equation}
Here we may and do assume that $m$ is independent of $i$. 
Replacing $P_{ij}(X)$ by its homogeneous degree $m-d$ part,
$Q_{ij}(X)$ by its homogeneous degree $m-d_-$ part,
and $R_i(X,T)$ by its homogeneous degree $m-1$ part, 
we may and do assume that the $P_{ij}(X)$, $Q_{ij}(X)$ and $R_i(X,T)$ are 
homogeneous polynomials with degree $m-d$, $m-d_-$ and $m-1$, 
respectively. 

Let $C'$ be the maximum of the absolute value of 
all the coefficients of $P_{ij}(X)$, $Q_{ij}(X)$ and $R_i(X,T)$ 
for $1 \leq i \leq N$ and $1 \leq j \leq N$.  
We set 
\begin{equation}
\label{eqn:C:def}
  C = \max\{C', 1\}.
\end{equation}

We fix real numbers $\varepsilon >0, \delta>0$ satisfying 
\begin{equation}
\label{eqn:e:d}
  \varepsilon \leq \frac{1}{C}, \quad 
  \delta \leq \frac{1}{C}, \quad
  (\varepsilon C)^d \leq \delta, \quad\text{and}\quad 
  (\varepsilon C)^{d_-} \leq \delta.
\end{equation}
Since $C \geq 1$, we have $0 < \varepsilon, \delta \leq 1$. 
For example, we can take 
\begin{equation}
\label{eqn:e:d:eg}
  \varepsilon = \frac{1}{C^{\min\{d, d_-\}}}  \quad\text{and}\quad
  \delta = \frac{1}{C^{\min\{d, d_-\}(\min\{d, d_-\}-1)}}.
\end{equation}

We define 
$N_{\delta, \varepsilon}^+$ and 
${V}_{\delta, \varepsilon}^+$ by  
\begin{align}
\label{eqn:n:v:+} 
  N_{\delta, \varepsilon}^+ 
  & := \left\{
  x \in \Aff^N(\Omega) \;\left\vert\;
  1 < \varepsilon \Vert x\Vert \quad\text{and}\quad
  \max\{\Vert f(x)\Vert, 1\} < \delta \max\{\Vert x\Vert^d, 1\}
  \right.\right\}, \\
\notag  
  {V}_{\delta, \varepsilon}^+ 
  & := \Aff^N(\Omega) \setminus N_{\delta, \varepsilon}^+. 
\end{align}

\begin{Remark}
We set 
\[
  \overline{N}_{\delta, \varepsilon}^+
  = \left\{
  (x:t) \in \PP^N(\Omega) \;\left\vert\;
  |t| < \varepsilon \Vert x\Vert \quad\text{and}\quad 
  \Vert(F(x,t),t^d)\Vert < \delta\Vert(x,t)\Vert^d \right.
  \right\}. 
\]
Then $N_{\delta, \varepsilon}^+ 
= \overline{N}_{\delta, \varepsilon}^+ \cap \Aff^N(\Omega)$.  
If $(x:t) \in I_+$, then $t=0$ and $F(x,t)=0$. 
Thus $|t|=0$ and $\Vert(F(x,t),t^d)\Vert=0$, so we have 
\[
  I_+ \subseteq \overline{N}_{\delta, \varepsilon}^+. 
\]
Intuitively, points in $N_{\delta, \varepsilon}^+$ are 
near to the hyperplane $\{(x:0) \in \PP^N(\Omega)\}$ at 
infinity (measured by $\varepsilon$), and also near to $I_+$ in 
``the direction of $x$'' (measured by $\delta$). 
\end{Remark}

\begin{Lemma}
\label{lemma:invariance:of:V}
$f({V}_{\delta, \varepsilon}^+) \subseteq {V}_{\delta, \varepsilon}^+$. 
\end{Lemma}
\Proof
Taking the complement, it suffices to show that 
\[
  f^{-1}(N_{\delta, \varepsilon}^+) \subseteq N_{\delta, \varepsilon}^+. 
\]
Suppose that $x=(x_1, \ldots, x_N) \in N_{\delta, \varepsilon}^+$. 
Without loss of generality, we may assume that $|x_1| = \Vert x\Vert$. 
We note that $f(x) = (F_1(x,1), \ldots, F_N(x,1))$ and 
$f^{-1}(x) = (G_1(x,1), \ldots, G_N(x,1))$. 
Since $\varepsilon \leq 1$, we have $\Vert x\Vert >1$. 
Then the definition of $N_{\delta, \varepsilon}^+$ gives  
\begin{align}
  \label{eqn:condition:1}
  \frac{1}{\varepsilon} & < \Vert x\Vert, \\
  \label{eqn:condition:2}
  \Vert f(x)\Vert & < \delta \Vert x\Vert^d. 
\end{align}
We need to show that $f^{-1}(x)\in N_{\delta, \varepsilon}^+$, 
which is equivalent to 
\begin{gather}
  \label{eqn:condition:3}
  1 < \varepsilon \Vert f^{-1}(x)\Vert, \\
  \label{eqn:condition:4}
  \max\{\Vert x\Vert, 1\} < \delta \max\{\Vert f^{-1}(x)\Vert^d, 1\}. 
  \end{gather}
First we show \eqref{eqn:condition:3}. 
To derive a contradiction, we assume that 
$\Vert f^{-1}(x)\Vert \leq \frac{1}{\varepsilon}$. 
Let $\lambda>0$ be any small number. 
We have 
{\allowdisplaybreaks
\begin{align*}
  & \left|\sum_{j=1}^N P_{1j}(x)F_j(x,1) 
  + \sum_{j=1}^N Q_{1j}(x)G_j(x,1) +  R_1(x,1)\right| \\
  & \qquad < 
  \max\left\{
   C\Vert x \Vert^{m-d}\cdot \delta \Vert x \Vert^d, \;
  (C+\lambda)\Vert x \Vert^{m-d_-}\frac{1}{\varepsilon}, \;
  (C+\lambda)\Vert x \Vert^{m-1}
  \right\} & \\
  & \qquad \leq 
  \max\left\{
  C \delta \Vert x \Vert^{m}, \;
  (C+\lambda)\Vert x \Vert^{m-d_-+1}, \;
  (C+\lambda)\Vert x \Vert^{m-1} 
  \right\} & \text{(from \eqref{eqn:condition:1})}\\
  & \qquad \leq 
  \max\left\{
   C \delta \Vert x \Vert^{m}, \;
  (C+\lambda)\Vert x \Vert^{m-1} 
  \right\} & \text{(since $d_- \geq 2$)}.
\end{align*}
}
Since $\lambda>0$ is arbitrary, 
\eqref{eqn:pqr} and the assumption that $|x_1| = \Vert x\Vert$ 
then gives either $\Vert x\Vert^m \leq C \Vert x \Vert^{m-1}$ 
or $\Vert x\Vert^m < C \delta \Vert x\Vert^m$. Equivalently, 
we have either $\Vert x\Vert \leq C$ or $1 < C \delta$. However, 
the former contradicts with \eqref{eqn:e:d} 
and \eqref{eqn:condition:1}, while the latter contradicts with 
\eqref{eqn:e:d}. Hence we get \eqref{eqn:condition:3}. 

Next we show \eqref{eqn:condition:4}. 
Since $\Vert x \Vert > \frac{1}{\varepsilon} \geq 1$ from 
\eqref{eqn:condition:1} and $\Vert f^{-1}(x) \Vert  > 
\frac{1}{\varepsilon} \geq 1$ from \eqref{eqn:condition:3}, 
the condition \eqref{eqn:condition:4} is equivalent to 
$\Vert x \Vert < \delta \Vert f^{-1}(x)\Vert^d$. To derive 
a contradiction, we assume the contrary, i.e., 
$\Vert x \Vert \geq \delta \Vert f^{-1}(x)\Vert^d$. 
Letting $\lambda >0$ be any small number, we have 
{\allowdisplaybreaks
\begin{align*}
  & \left|\sum_{j=1}^N P_{1j}(x)F_j(x,1) 
  + \sum_{j=1}^N Q_{1j}(x)G_j(x,1) +  R_1(x,1)\right| \\
  & \qquad < 
  \max\left\{
   C\Vert x \Vert^{m-d}\cdot \delta \Vert x \Vert^d, \;
  (C+\lambda)\Vert x \Vert^{m-d_-}\cdot 
  \left(\frac{1}{\delta}\right)^{\frac{1}{d}} \Vert x \Vert^{\frac{1}{d}}, \;
  (C+\lambda)\Vert x \Vert^{m-1}
  \right\}  \\
  & \qquad \leq 
  \max\left\{
  C \delta \Vert x \Vert^{m}, \;
  (C+\lambda)\left(\frac{1}{\delta}\right)^{\frac{1}{d}}
  \Vert x \Vert^{m-d_-+\frac{1}{d}}, \;
  (C+\lambda)\Vert x \Vert^{m-1} 
  \right\} \\
  & \qquad \leq 
  \max\left\{
   C \delta \Vert x \Vert^{m}, \;
  (C+\lambda)\left(\frac{1}{\delta}\right)^{\frac{1}{d}}\Vert x \Vert^{m-1} 
  \right\} \qquad\qquad\qquad\qquad \text{(since $d_- - \frac{1}{d} \geq 1$)}.  
\end{align*}
}
Since $\lambda>0$ is arbitrary, 
\eqref{eqn:pqr} and the assumption that $|x_1| = \Vert x\Vert$ 
gives this time 
\[
  \text{either}\quad 
  \Vert x\Vert \leq \left(\frac{1}{\delta}\right)^{\frac{1}{d}}C\quad 
  \text{or}\quad 
  1 < C \delta. 
\]
However, the former contradicts with \eqref{eqn:e:d} 
and \eqref{eqn:condition:1}, while the latter contradicts with 
\eqref{eqn:e:d}. Hence we get \eqref{eqn:condition:4}, which completes 
the proof. 
\QED

\begin{Lemma}
\label{lemma:before:G}
Set 
$C_{\delta, \varepsilon}^+ := \min\{\delta, \varepsilon^d\}$. 
Then 
\[
  \max\{\Vert f(x) \Vert, 1\} 
  \geq C_{\delta, \varepsilon}^+ \cdot \max\{\Vert x \Vert^d, 1\}
  \qquad\text{for all $x \in {V}_{\delta, \varepsilon}^+$.}
\]
\end{Lemma}

\Proof
For $x \in {V}_{\delta, \varepsilon}^+$,  
the definition of ${V}_{\delta, \varepsilon}^+$ gives 
\[
  \text{either}\quad 
  \Vert x\Vert \leq \frac{1}{\varepsilon}\quad
  \text{or}\quad 
  \max\{\Vert f(x)\Vert, 1\} \geq  \delta\max\{\Vert x\Vert^d, 1\}. 
\]
If the latter holds, then we get the assertion since 
$\delta \geq C_{\delta, \varepsilon}^+$. If the former holds, 
then $C_{\delta, \varepsilon}^+ \Vert x\Vert^d \leq 1$. 
Noting that $C_{\delta, \varepsilon}^+ \leq 1$, we get 
$\max\{\Vert f(x) \Vert, 1\} 
  \geq 1 
  \geq C_{\delta, \varepsilon}^+ \cdot \max\{\Vert x \Vert^d, 1\}$. 
\QED

\begin{Lemma}
\label{lemma:estimate:G}
Set $c_{\delta, \varepsilon}^+ := 
\frac{1}{d-1} \log C_{\delta, \varepsilon}^+$. 
Then 
\[
  G_f(x)
  \geq \log^+\Vert x \Vert + c_{\delta, \varepsilon}^+
  \qquad\text{for all $x \in {V}_{\delta, \varepsilon}^+$.}
\]
\end{Lemma}

\Proof
Suppose $x \in {V}_{\delta, \varepsilon}^+$. 
It follows from Lemma~\ref{lemma:invariance:of:V} that  
$f^n(x) \in {V}_{\delta, \varepsilon}^+$ for all $n \geq 1$. 
Then Lemma~\ref{lemma:before:G} gives 
\[
  \log^+\Vert f^{n}(x)\Vert \geq d \log^+\Vert f^{n-1}(x)\Vert 
  + \log C_{\delta, \varepsilon}^+. 
\]
The usual telescoping argument tells us that 
{\allowdisplaybreaks
\begin{align*}
  G_f(x) 
  & = \lim_{n \to +\infty} \frac{1}{d^n} \log^+\Vert f^n(x)\Vert \\
  & = \log^+\Vert x \Vert+ 
    \sum_{n=1}^{\infty} \frac{1}{d^n} 
    \left(\log^+\Vert f^{n}(x)\Vert -  d \log^+\Vert f^{n-1}(x)\Vert\right) \\
  & \geq  \log^+\Vert x \Vert + c_{\delta, \varepsilon}^+. 
\end{align*} 
}
This completes the proof. 
\QED

With $f^{-1}$ in place of $f$, 
we define $N_{\delta, \varepsilon}^-$ and 
${V}_{\delta, \varepsilon}^-$ by  
\begin{align}
\label{eqn:n:v:-} 
N_{\delta, \varepsilon}^-
& := \left\{
x \in \Aff^N(\Omega) \;\left\vert\;
1 < \varepsilon \Vert x\Vert \quad\text{and}\quad 
\max\{\Vert f^{-1}(x)\Vert, 1\} < \delta\max\{\Vert x\Vert^{d_-} , 1\}
\right.
\right\}
\\
\notag
{V}_{\delta, \varepsilon}^- 
& := \Aff^N(\Omega) \setminus N_{\delta, \varepsilon}^-. 
\end{align}
Then, setting $c_{\delta, \varepsilon}^- := 
\frac{1}{d_--1} \log\min\{\delta, \epsilon^{d_-}\}$, we have 
\begin{equation}
\label{eqn:estimate:G-}
  G_{f^{-1}}(x)
  \geq \log^+\Vert x \Vert + c_{\delta, \varepsilon}^-
  \qquad\text{for all $x \in {V}_{\delta, \varepsilon}^-$.}
\end{equation}

\begin{Lemma}
\label{lemma:V:V-}
${V}_{\delta, \varepsilon}^+ \cup {V}_{\delta, \varepsilon}^- 
= \Aff^N(\Omega)$. 
\end{Lemma}

\Proof
Taking the complement, it suffices to show that 
$N_{\delta, \varepsilon}^+ \cap N_{\delta, \varepsilon}^- = \emptyset$.
To derive a contradiction, we assume that there is an 
$x \in N_{\delta, \varepsilon}^+ \cap N_{\delta, \varepsilon}^-$. 
Then we have 
{\allowdisplaybreaks
\begin{align}
\label{eqn:V:V-:1}
\Vert x \Vert & > \frac{1}{\varepsilon} \\
\label{eqn:V:V-:2}
\max\{\Vert f(x) \Vert, 1\} & <  \delta \max\{\Vert x\Vert^d, 1\} \\
\label{eqn:V:V-:3}
\max\{\Vert f^{-1}(x) \Vert, 1\} & <  \delta \max\{\Vert x\Vert^{d_-}, 1\} 
\end{align}
}
Without loss of generality, we assume that $|x_1| = \Vert x\Vert$. 
Let $\lambda >0$ be any small number. By 
\eqref{eqn:V:V-:1}--\eqref{eqn:V:V-:3}, 
we have 
{\allowdisplaybreaks
\begin{align*}
  & \left|\sum_{j=1}^N P_{1j}(x)F_j(x,1) 
  + \sum_{j=1}^N Q_{1j}(x)G_j(x,1) +  R_1(x,1)\right| \\
  & \qquad < 
  \max\left\{
   C\Vert x \Vert^{m-d}\cdot \delta \Vert x \Vert^d, \;
   C\Vert x \Vert^{m-d_-}\cdot \delta \Vert x \Vert^{d_-}, \;
  (C+\lambda)\Vert x \Vert^{m-1}
  \right\}  \\
  & \qquad \leq 
  \max\left\{
  C \delta \Vert x \Vert^{m}, \;
  (C+\lambda)\Vert x \Vert^{m-1}
   \right\}. 
\end{align*}
}
Since $\lambda$ is arbitrary, 
it follows from \eqref{eqn:pqr} that 
either $\Vert x\Vert^m < C\delta \Vert x \Vert^m$ or 
$\Vert x\Vert^m \leq C\Vert x \Vert^{m-1}$. Hence we get 
\[
  \text{either}\quad 
  1 < C\delta \quad 
  \text{or}\quad 
  \Vert x\Vert \leq C. 
\]
However, the former contradicts with \eqref{eqn:e:d}, 
while the latter contradicts with 
\eqref{eqn:e:d} and \eqref{eqn:V:V-:1}. 
Thus we have $N_{\delta, \varepsilon}^+ 
\cap N_{\delta, \varepsilon}^- = \emptyset$.
\QED

\medskip
{\sl Proof of Theorem~\ref{thm:local:II}}. 
Let $\delta$ and $\varepsilon$ be constants satisfying 
\eqref{eqn:e:d}. Then Theorem~\ref{thm:local:II} holds 
with ${V}^\pm = {V}_{\delta, \varepsilon}^\pm$ and 
$c^\pm = c_{\delta, \varepsilon}^\pm$. 
Indeed, the condition (i) follows from Lemma~\ref{lemma:estimate:G}, 
and the condition (ii) from \eqref{eqn:estimate:G-}, while 
the condition (iii) follows from Lemma~\ref{lemma:V:V-}. 
\QED

\medskip
\section{Non-archimedean Green functions and the set of 
escaping points}
\label{sec:local:theory:III}
\setcounter{equation}{0}
In this section, we continue to study some basic properties of 
regular polynomial automorphisms defined over $\Omega$. 
We keep the notation and the assumption of \S\ref{sec:local:theory:II}. 
In particular, $f: \Aff^N\to\Aff^N$ denotes a regular 
polynomial automorphism of degree $d \geq 2$ defined over 
$\Omega$. 

In analogy with the field of complex numbers, 
we define the set $W^+$ of escaping points and 
the set $\Kcal^+$ of non-escaping points by  
\begin{align*}
  W^+ & := \{x \in \Aff^N(\Omega) \mid \Vert f^n(x)\Vert \to +\infty 
      \;(n\to +\infty) \}, \\
  \Kcal^+ & :=  \{x \in \Aff^N(\Omega) \mid  
       \text{$\{f^n(x)\}_{n=0}^{+\infty}$ is bounded with respect to 
       $\Vert\cdot\Vert$}\} 
\end{align*}
Then the following theorem holds, which is a non-archimedean 
version of the results of \cite[\S\S2--3]{BS} and \cite[\S2]{Sib}. 

\begin{Theorem}
\label{thm:G:K}
Let $f: \Aff^N\to\Aff^N$ be a regular polynomial automorphism 
over $\Omega$, and let $G_f$ be the Green function for $f$. 
\begin{enumerate}
\item[(1)]
The set $\Kcal^+$ is exactly the set of points where $G_f$ vanish\textup{:} 
\[
  \Kcal^+ = \{x \in \Aff^N(\Omega) \mid G_f(x) = 0\}.
\] 
\item[(2)]
$\Aff^N(\Omega) = W^+ \amalg \Kcal^+$ 
\textup{(}disjoint union\textup{)}.  
\end{enumerate}
\end{Theorem}

To prove Theorem~\ref{thm:G:K}, 
we need the following two lemmas. Recall that $\delta$ and $\varepsilon$ 
are fixed constants satisfying \eqref{eqn:e:d}. 

\begin{Lemma}
\label{lemma:for:G:K:1}
For any $x \in N_{\delta, \frac{\varepsilon}{2}}^+$, 
one has $\Vert x \Vert \leq \frac{1}{2} \Vert f^{-1}(x)\Vert$. 
\end{Lemma}

\Proof
It follows from $x \in N_{\delta, \frac{\varepsilon}{2}}^+$ that 
\begin{equation}
\label{eqn:for:G:K:1}
  \Vert x\Vert > \frac{2}{\varepsilon}
  \quad\text{and}\quad
  \Vert f(x)\Vert < \delta \Vert x\Vert^d. 
\end{equation}
To derive a contradiction, 
we assume that $\Vert x \Vert > \frac{1}{2} \Vert f^{-1}(x)\Vert$. 
Without loss of generality, 
we assume that $|x_1| = \Vert x \Vert$. 
Then  (we take $\lambda = C$ here)
{\allowdisplaybreaks
\begin{align*}
  & \left|\sum_{j=1}^N P_{1j}(x)F_j(x,1) 
  + \sum_{j=1}^N Q_{1j}(x)G_j(x,1) +  R_1(x,1)\right| \\
  & \qquad < 
  \max\left\{
   C\Vert x \Vert^{m-d}\cdot \delta \Vert x \Vert^d, \;
   C\Vert x \Vert^{m-d_-}\cdot 2 \Vert x \Vert, \;
   2C \Vert x \Vert^{m-1}
  \right\}  \\
  & \qquad \leq 
  \max\left\{
  C \delta \Vert x \Vert^{m}, \;
  2 C \Vert x \Vert^{m-1}
   \right\}. 
\end{align*}
}
Using \eqref{eqn:pqr}, we get 
\[
  \text{either}\quad
  1 < C\delta
  \quad\text{or}\quad
  \Vert x\Vert < 2C. 
\]
However, the former contradicts with \eqref{eqn:e:d}. If the latter holds, 
then Eqn.~\eqref{eqn:for:G:K:1} implies $1< C\varepsilon$, 
contradicting with \eqref{eqn:e:d}. 
This completes the proof. 
\QED

\begin{Lemma}
\label{lemma:for:G:K:2}
For any $x \in \Aff^N(\Omega)$, 
one has $f^n(x) \in 
V_{\delta, \frac{\varepsilon}{2}}^+$ for 
all sufficiently large $n$. 
\end{Lemma}

\Proof
Note that $\frac{\varepsilon}{2}, \delta$ satisfies 
\eqref{eqn:e:d} with $\frac{\varepsilon}{2}$ 
in place of $\varepsilon$. Thus, 
if $x \in V_{\delta, \frac{\varepsilon}{2}}^+$, 
then Lemma~\ref{lemma:invariance:of:V} 
gives $f^n(x) \in V_{\delta, \frac{\varepsilon}{2}}^+$ for all $n \geq 0$. 

Suppose now that $x \in N_{\delta, \frac{\varepsilon}{2}}^+$. 
We take a positive integer $n_0$ so that 
$\Vert x \Vert \leq \frac{2^{n_0+1}}{\varepsilon}$. 
We claim that $f^{n_0}(x) \in V_{\delta, \frac{\varepsilon}{2}}^+$. 
Indeed, if we assume the contrary, then 
Lemma~\ref{lemma:for:G:K:1}, applied to 
$x, \ldots, f^{n_0}(x) \in N_{\delta, \frac{\varepsilon}{2}}^+$,  
gives 
\[
  \frac{2}{\varepsilon} 
  <  \Vert f^{n_0}(x)\Vert 
  \leq \frac{1}{2} \Vert f^{n_0-1}(x)\Vert 
  \leq\cdots\leq \frac{1}{2^{n_0}} \Vert x\Vert,     
\]
which contradicts with our choice of $n_0$. 
Thus $f^n(x) \in V_{\delta, \frac{\varepsilon}{2}}^+$ for 
all $n \geq n_0$.
\QED

{\sl Proof of Theorem~\ref{thm:G:K}.}\quad
(1) From Definition~\ref{def:Green:function}, we get 
$\Kcal^+ \subseteq \{x\in\Aff^N(\Omega) \mid G_f(x) =0\}$. To show the other 
inclusion, we assume that $G_f(x)=0$. Then $G_f(f^n(x)) = d^n G_f(x) = 0$ 
for all $n \geq 0$. By Lemma~\ref{lemma:for:G:K:2}, 
we take $n_0$ such that $f^{n_0}(x) \in V_{\delta, \frac{\varepsilon}{2}}^+$. 
It follows from Lemma~\ref{lemma:invariance:of:V} and 
Lemma~\ref{lemma:estimate:G} (applied to $\frac{\varepsilon}{2}$ in place of $\varepsilon$) that 
\[
  G_f(f^n(x)) \geq \log^+\Vert f^n(x)\Vert + 
  c_{\delta, \frac{\varepsilon}{2}}^+ 
\]
for all $n \geq n_0$. Combined with 
$G_f(f^n(x)) = 0$, we see that 
$\Vert f^n(x) \Vert \leq 
\exp\left(- c_{\delta, \frac{\varepsilon}{2}}^+\right)$ 
for all $n \geq n_0$. 
Thus $\{x\in\Aff^N(\Omega) \mid G_f(x) =0\} \subseteq \Kcal^+$.

\medskip
(2) If $x \not\in \Kcal^+$, then $G_f(x) >0$ by (1). 
Definition~\ref{def:Green:function} then gives 
$\Vert f^n(x)\Vert \to +\infty$ as $n \to +\infty$.  
\QED

With $f^{-1}$ in place of $f$, we put 
\begin{align*}
  W^- & := \{x \in \Aff^N(\Omega) \mid \Vert f^{-n}(x)\Vert \to +\infty 
      \;(n\to +\infty) \}, \\
  \Kcal^- & :=  \{x \in \Aff^N(\Omega) \mid  
       \text{$\{f^{-n}(x)\}_{n=0}^{+\infty}$ is bounded with respect to 
       $\Vert\cdot\Vert$}\} 
\end{align*}
Then we have $\Aff^N(\Omega) = W^- \amalg \Kcal^-$ as in Theorem~\ref{thm:G:K}.

\smallskip
In the rest of this section, we give filtrations 
of $\Aff^N$ relative to $f$ over non-archimedean fields, as in
Bedford--Smillie \cite[\S2.2]{BS} and \cite[\S3]{SW} over $\CC$. 

We set 
\begin{align*}
  B_{\delta, \varepsilon} 
  & = \left\{x \in \Aff^N(\Omega) \;\left\vert\; 
  \Vert x \Vert \leq \frac{1}{\varepsilon} \right.\right\}, \\
  U^+_{\delta, \varepsilon} 
  & = \left\{x \in \Aff^N(\Omega) \;\left\vert\; 
  \text{$\Vert x \Vert > \frac{1}{\varepsilon}$ 
  and $\Vert f(x) \Vert \geq \delta \Vert x \Vert^d$} \right.\right\}, \\
  U^-_{\delta, \varepsilon} 
  & = \left\{x \in \Aff^N(\Omega) \;\left\vert\; 
  \text{$\Vert x \Vert > \frac{1}{\varepsilon}$ 
  and $\Vert f(x) \Vert < \delta \Vert x \Vert^d$} 
  \right.\right\}, 
\end{align*}
where $\delta$ and $\varepsilon$ are constants satisfying 
\eqref{eqn:e:d}. 

Since $\varepsilon \leq 1$ 
and $\frac{\delta}{\varepsilon^d} \geq C^d \geq 1$ 
by \eqref{eqn:e:d}, we have 
\[
  U^+_{\delta, \varepsilon} 
  = \left\{x \in \Aff^N(\Omega) \;\left\vert\; 
  \text{$\Vert x \Vert > \frac{1}{\varepsilon}$ 
  and $\max\{\Vert f(x) \Vert, 1\} \geq \delta 
  \max\{\Vert x \Vert, 1\}^d$} \right.\right\}, 
\]
so that $B_{\delta, \varepsilon} \amalg U^+_{\delta, \varepsilon} 
= V^+_{\delta, \varepsilon}$.

\begin{Proposition}
\label{prop:filtration}
We assume that $\varepsilon$ and $\delta$ satisfy 
\begin{equation}
\label{eqn:e:d:2}
  \varepsilon^{d-1} \leq  \delta 
  \quad\text{and}\quad
 \varepsilon^{d_--1} \leq \delta 
\end{equation}
in addition to $\eqref{eqn:e:d}$ \textup{(}For example, if we take $\varepsilon$ and 
$\delta$ as \eqref{eqn:e:d:eg}, then they also satisfy  \eqref{eqn:e:d:2}\textup{)}. 
Then we have the following. 
\begin{enumerate}
\item[(1)] 
$\Aff^N(\Omega) = B_{\delta, \varepsilon}\amalg 
U^+_{\delta, \varepsilon}\amalg  U^-_{\delta, \varepsilon}$ 
\textup{(}disjoint union\textup{)}. 
\item[(2)]
$f(U^+_{\delta, \varepsilon}) \subseteq U^+_{\delta, \varepsilon}$ 
and $f(B_{\delta, \varepsilon} \amalg U^+_{\delta, \varepsilon}) 
\subseteq B_{\delta, \varepsilon} \amalg U^+_{\delta, \varepsilon}$.  
\item[(3)]
$f^{-1}(U^-_{\delta, \varepsilon}) \subseteq U^-_{\delta, \varepsilon}$ 
and $f^{-1}(B_{\delta, \varepsilon} \amalg U^-_{\delta, \varepsilon}) 
\subseteq B_{\delta, \varepsilon} \amalg U^-_{\delta, \varepsilon}$.  
\end{enumerate}
\end{Proposition}

\Proof
(1) This is obvious from the definition. 

\smallskip
(2) Since $B_{\delta, \varepsilon} \amalg U^+_{\delta, \varepsilon} 
= V^+_{\delta, \varepsilon}$, we have 
$f(B_{\delta, \varepsilon} \amalg U^+_{\delta, \varepsilon}) 
\subseteq B_{\delta, \varepsilon} \amalg U^+_{\delta, \varepsilon}$ 
by Lemma~\ref{lemma:invariance:of:V}. 

Suppose that $x \in U^+_{\delta, \varepsilon}$. Then 
\begin{equation}
\label{eqn:x:U+}
  \Vert f(x) \Vert 
  \geq \delta \Vert x \Vert^d > \frac{\delta}{\varepsilon^d} 
  \geq \frac{1}{\varepsilon},  
\end{equation}
where we have used \eqref{eqn:e:d:2} in the last inequality. 
Also, since $x \in U^+_{\delta, \varepsilon} 
\subseteq V^+_{\delta, \varepsilon}$, we have  
$f(x) \in V^+_{\delta, \varepsilon}$ by 
Lemma~\ref{lemma:invariance:of:V}. Since $f(x) 
\not\in B_{\delta, \varepsilon}$ by \eqref{eqn:x:U+}, 
we get $f(x) \in V^+_{\delta, \varepsilon}\setminus 
B_{\delta, \varepsilon} = U^+_{\delta, \varepsilon}$. 
Hence $f(U^+_{\delta, \varepsilon}) \subseteq U^+_{\delta, \varepsilon}$. 

\smallskip
(3) We put 
\begin{align}
\label{eqn:x:U:tilde}
  \widetilde{U}^-_{\delta, \varepsilon}
  & = \left\{x \in \Aff^N(\Omega) \;\left\vert\; 
  \text{$\Vert x \Vert > \frac{1}{\varepsilon}$ 
  and $\Vert f^{-1}(x) \Vert \geq \delta \Vert x \Vert^{d_-}$} 
  \right.\right\} \\
\notag
  & = \left\{x \in \Aff^N(\Omega) \;\left\vert\; 
  \text{$\Vert x \Vert > \frac{1}{\varepsilon}$ 
  and $\max\{\Vert f^{-1}(x) \Vert, 1\} \geq \delta 
  \max\{\Vert x \Vert, 1\}^{d_-}$} \right.\right\},  
\end{align}
where the second equality follows from $\frac{\delta}{\varepsilon^d_-} \geq C^{d_-} \geq 1$ 
by \eqref{eqn:e:d}. 
Then as in (2), we have 
$f^{-1}(\widetilde{U}^-_{\delta, \varepsilon}) 
\subseteq \widetilde{U}^-_{\delta, \varepsilon}$. 
Since $B_{\delta, \varepsilon} \amalg \widetilde{U}^-_{\delta, \varepsilon} 
= V^-_{\delta, \varepsilon}$, Lemma~\ref{lemma:V:V-} implies 
that $U^-_{\delta, \varepsilon} \subseteq 
\widetilde{U}^-_{\delta, \varepsilon}$. 

Suppose that $x \in U^-_{\delta, \varepsilon}$. Then 
\[
  f^{-1}(x) \in f^{-1}(U^-_{\delta, \varepsilon}) 
  \subseteq  f^{-1}(\widetilde{U}^-_{\delta, \varepsilon}) 
  \subseteq \widetilde{U}^-_{\delta, \varepsilon}. 
\]
In particular, $\Vert f^{-1}(x) \Vert > \frac{1}{\varepsilon}$, 
so that $f^{-1}(x) \not\in B_{\delta, \varepsilon}$. 
On the other hand, since $x \not\in U^+_{\delta, \varepsilon}$ 
and $f(U^+_{\delta, \varepsilon}) \subseteq U^+_{\delta, \varepsilon}$, 
we get $f^{-1}(x) \not\in U^+_{\delta, \varepsilon}$. 
It follows that $f^{-1}(x) \in U^-_{\delta, \varepsilon} 
= \Aff^N(\Omega) \setminus (B_{\delta, \varepsilon} 
\amalg U^+_{\delta, \varepsilon})$. We conclude that 
$f^{-1}(U^-_{\delta, \varepsilon}) \subseteq U^-_{\delta, \varepsilon}$. 

Next we show $f^{-1}(B_{\delta, \varepsilon} 
\amalg U^-_{\delta, \varepsilon}) 
\subseteq B_{\delta, \varepsilon} \amalg U^-_{\delta, \varepsilon}$. 
Since $U^+_{\delta, \varepsilon} = \Aff^N(\Omega) \setminus (B_{\delta, \varepsilon} 
\amalg U^-_{\delta, \varepsilon})$, it suffices to show 
that $f^{-1}(U^+_{\delta, \varepsilon}) 
\supseteq U^+_{\delta, \varepsilon}$. This inclusion is obvious from 
$f(U^+_{\delta, \varepsilon}) \subseteq U^+_{\delta, \varepsilon}$. 
\QED

\begin{Proposition}
\label{prop:filtration:2}
We assume that $\varepsilon$ and $\delta$ satisfy 
\begin{equation}
\label{eqn:e:d:2:s}
  \varepsilon^{d-1} <  \delta 
  \quad\text{and}\quad
  \varepsilon^{d_--1} < \delta 
\end{equation}
in addition to $\eqref{eqn:e:d}$. Then we have the following. 
\begin{enumerate}
\item[(1)]
$\bigcup_{n=0}^{+\infty} f^{-n}(U^+_{\delta, \varepsilon})
= W^+$. 
\item[(2)]
$\bigcup_{n=0}^{+\infty} f^{n}(U^-_{\delta, \varepsilon})
= W^-$. 
\end{enumerate}
\end{Proposition}

\Proof
(1) We set $r := \frac{\delta}{\varepsilon^{d-1}} > 1$. 
We first show that $U^+_{\delta, \varepsilon} \subseteq W^+$. 
Indeed, if $x \in U^+_{\delta, \varepsilon}$, then 
\[ 
  \Vert f(x) \Vert \geq \delta \Vert x \Vert^d  
  > \frac{\delta}{\varepsilon^{d-1}} \frac{1}{\varepsilon} 
  = r \frac{1}{\varepsilon}. 
\]
Since $f(U^+_{\delta, \varepsilon}) 
\subseteq U^+_{\delta, \varepsilon}$, we inductively get 
$\Vert f^n(x) \Vert > r^{\frac{d^n-1}{d-1}} \frac{1}{\varepsilon}$ 
for all $n \geq 0$. Hence $x \in W^+$. This completes the 
proof of $U^+_{\delta, \varepsilon} \subseteq W^+$. 
 Since $f^{-1}(W^+) = W^+$, we get $f^{-n}(U^+_{\delta, \varepsilon}) 
\subseteq W^+$ for all $n \geq 0$, so that  
$\bigcup_{n=0}^{+\infty} f^{-n}(U^+_{\delta, \varepsilon}) \subseteq W^+$. 

To show the other inclusion $\bigcup_{n=0}^{+\infty} 
f^{-n}(U^+_{\delta, \varepsilon}) \supseteq W^+$, suppose that 
$x \not\in \bigcup_{n=0}^{+\infty} f^{-n}(U^+_{\delta, \varepsilon})$. 
We need to show that $x \in \Kcal^+$. Since $f^n(x) \not\in 
U^+_{\delta, \varepsilon}$, we have 
either $f^n(x) \in B_{\delta, \varepsilon}$ or 
$f^n(x) \in U^-_{\delta, \varepsilon}$. 

{\bf Case~1:}\quad Suppose that there is an $n_0 \geq 0$ such that $f^{n_0}(x) 
\in B_{\delta, \varepsilon}$. Then $f^{n_0+1}(x) \in B_{\delta, \varepsilon} 
\amalg U^+_{\delta, \varepsilon}$ by Proposition~\ref{prop:filtration}(2). 
Since $f^{n_0+1}(x) \not\in 
U^+_{\delta, \varepsilon}$, we obtain $f^{n_0+1}(x) 
\in B_{\delta, \varepsilon}$. Inductively, $f^{n}(x) 
\in B_{\delta, \varepsilon}$ for all $n \geq n_0$, so we conclude that 
$x \in \Kcal^+$. 

{\bf Case~2:}\quad Suppose that $f^{n}(x) \in U^-_{\delta, \varepsilon}$ for 
all $n \geq 0$. By Lemma~\ref{lemma:for:G:K:2}, there is an 
$n_0 \geq 0$ such that $f^n(x) \in V_{\delta, \frac{\varepsilon}{2}}^+$ 
for all $n \geq n_0$. Then for all $n \geq n_0$, we have 
\[
  f^n(x) \in V_{\delta, \frac{\varepsilon}{2}}^+ 
  \cap U^-_{\delta, \varepsilon} \subseteq 
  \left\{y \in \Aff^N(\Omega) \;\left\vert\; \frac{1}{\varepsilon} 
  < \Vert y \Vert \leq \frac{2}{\varepsilon}\right.\right\}. 
\]
Hence $x \in \Kcal^+$. 

In both cases we have $x \in \Kcal^+$, so we get $\bigcup_{n=0}^{+\infty} 
f^{-n}(U^+_{\delta, \varepsilon}) \supseteq W^+$. 

\smallskip
(2) 
Let $\widetilde{U}^-_{\delta, \varepsilon}$ be the set defined by \eqref{eqn:x:U:tilde}. 
Then the argument in (1) gives $\bigcup_{n=0}^{+\infty} 
f^{n}(\widetilde{U}^-_{\delta, \varepsilon}) = W^-$, and so 
$\bigcup_{n=0}^{+\infty} 
f^{n}(U^-_{\delta, \varepsilon}) \subseteq W^-$. 
To show the other inclusion, suppose that 
$x \not\in \bigcup_{n=0}^{+\infty} f^{n}(U^-_{\delta, \varepsilon})$. 
Then we have 
either $f^{-n}(x) \in B_{\delta, \varepsilon}$ or 
$f^{-n}(x) \in U^+_{\delta, \varepsilon}$. 

{\bf Case~1:}\quad 
If there is an $n_0 \geq 0$ such that $f^{-n_0}(x) 
\in B_{\delta, \varepsilon}$, then the argument of Case~1 of (1), 
together with Proposition~\ref{prop:filtration}(3), gives 
$f^{-n}(x) \in B_{\delta, \varepsilon}$ for all $n \geq n_0$. 

{\bf Case~2:}\quad 
Suppose that $f^{-n}(x) \in U^+_{\delta, \varepsilon}$ for 
all $n \geq 0$.

then the argument of Case~2 of (4), 
together with Lemma~\ref{lemma:for:G:K:2} with $f^{-1}$ in place 
of $f$, gives $\frac{1}{\varepsilon} 
  < \Vert x \Vert < \frac{2}{\varepsilon}$ for sufficiently 
large $n$. 

In both cases, we get $x \in \Kcal^-$. 
Hence $\bigcup_{n=0}^{+\infty} 
f^{n}(U^-_{\delta, \varepsilon}) \supseteq W^-$. 
\QED

\begin{Remark}
If we take 
\[
  0 < \varepsilon < \frac{1}{C^{\min\{d, d_-\}}}
  \quad\text{and}\quad
  \delta = \frac{1}{C^{\min\{d, d_-\}(\min\{d, d_-\}-1)}}, 
\] 
then they satisfy both \eqref{eqn:e:d} and \eqref{eqn:e:d:2:s}. 
\end{Remark}

\medskip
\section{Regular automorphisms having good reduction}
\label{sec:good:reduction}
\setcounter{equation}{0}
In \cite{MS}, Morton and Silverman introduced 
the notion of having good reduction for endomorphisms of $\PP^1$ over 
$\Omega$, which has been useful in studying 
endomorphisms of $\PP^1$ over a global field. For 
endomorphisms of $\PP^N$ having good reduction, see, for example 
\cite[Remark~12]{KS} and \cite{KS2}. 
In this section, we introduce the notion of 
having good reduction for regular polynomial automorphisms 
of $\Aff^N$ over $\Omega$. This notion will be useful 
in studying regular polynomial automorphisms over a global field 
in \S\S\ref{sec:global:theory}--\ref{sec:applications}. 

As in \S\ref{sec:local:theory}, $R$ denotes the ring of integers of $\Omega$. 
Let $M$ be the maximal ideal of $R$, and $\widetilde{k} := R/M$ 
the residue field. Note that $\widetilde{k}$ is algebraically closed 
since $\Omega$ is algebraically closed. 

\begin{Definition}[good reduction]
\label{def:good:reduction}
Let $f=(f_1, \ldots, f_N): \Aff^N \to \Aff^N$ be a regular polynomial 
automorphism over an algebraically closed field $\Omega$ with non-trivial non-archimedean 
absolute value, and let $f^{-1}=(g_1, \ldots, g_N): 
\Aff^N \to \Aff^N$ denote its inverse. We write $d$ and $d_-$ for 
the degrees of $f$ and $f^{-1}$, respectively.  
We say that $f$ {\em has good reduction} if 
the following three conditions are satisfied. 
\begin{enumerate}
\item[(i)]
$f$ extends to the polynomial automorphism $f: \Aff^N_R \to \Aff^N_R$ 
over $R$, meaning that both $f_1(X), \ldots, f_N(X)$ and 
$g_1(X), \ldots, g_N(X)$ are in $R[X_1, \ldots, X_N]$. 
\item[(ii)]
Let $\widetilde{f}=(\widetilde{f_1}, \ldots, \widetilde{f_N}): 
\Aff^N_{\widetilde{k}} \to \Aff^N_{\widetilde{k}}$ 
and $\widetilde{f^{-1}}=(\widetilde{g_1}, \ldots, \widetilde{g_N}): 
\Aff^N_{\widetilde{k}} \to \Aff^N_{\widetilde{k}}$ be the induced 
polynomial automorphisms over $\widetilde{k}$. 
Then the degrees of $\widetilde{f}$ and $\widetilde{f^{-1}}$ are 
equal to $d$ and $d_-$, respectively. 
\item[(iii)]
$\widetilde{f}$ is regular (cf. Remark~\ref{rmk:regular:autom:def}). 
\end{enumerate}
\end{Definition}

We give some equivalent conditions for regular polynomial automorphisms $f$ 
to have good reduction. As in \S\ref{sec:local:theory}, 
let $F_i(X,T)$ and $G_j(X,T)$ be 
the homogenization of $f_i(X)$ and $g_j(X)$. If $F_i(X,T)$ and $G_j(X,T)$ 
are defined over $R$, we denote by $\widetilde{F_i}(X,T)$ and 
$\widetilde{G_j}(X,T)$ their reductions to $\widetilde{k}$. 
Let $\rho: R \to \widetilde{k}$ be the natural map. 

\begin{Proposition}
\label{prop:good:reduction:rephrasing}
Let $f$ be a regular polynomial automorphism 
of $\Aff^N$ over $\Omega$. 
Assume that $f$ satisfies the conditions \textup{(i)} and \textup{(ii)} of 
Definition~\textup{\ref{def:good:reduction}}. Then the followings are equivalent. 
\begin{enumerate}
\item[(1)]
$f$ has good reduction, i.e., $f$ also satisfies the condition \textup{(iii)} 
of Definition~\textup{\ref{def:good:reduction}}. 
\item[(2)]
As ideals in $R[X_1, \ldots, X_N, T]$, one has 
\[
  \left(X_1, \ldots, X_N, T\right)^k \subseteq 
  \left(F_1(X,T), \ldots, F_N(X,T), G_1(X,T), \ldots, G_N(X,T), T\right)
\]
for some integer $k \geq 1$. 
\item[(3)]
As ideals in $R[X_1, \ldots, X_N]$, one has 
\[
  \left(X_1, \ldots, X_N\right)^\ell \subseteq 
  \left(F_1(X,0), \ldots, F_N(X,0), G_1(X,0), \ldots, G_N(X,0)\right)
\]
for some integer $\ell \geq 1$. 
\end{enumerate}
\end{Proposition}

\Proof
(1) $\Longrightarrow$ (3): 
It suffices to show that 
\begin{equation}
\label{eqn:a:to:c}
  \left(X_1, \ldots, X_N\right)^\ell \subseteq 
  \left(F_1(X,0)^{d_-}, \ldots, F_N(X,0)^{d_-}, 
  G_1(X,0)^{d}, \ldots, G_N(X,0)^{d}\right)
\end{equation}
for some $\ell \geq 1$. We set  
\[
  I = \left\{r \in R \;\left\vert\; 
  \begin{gathered}
  \text{There is an $\ell \geq 1$ such that} \\
  r \left(X_1, \ldots, X_N\right)^\ell \subseteq 
  \left(F_1(X,0)^{d_-}, \ldots, F_N(X,0)^{d_-}, 
  G_1(X,0)^{d}, \ldots, G_N(X,0)^{d}\right)
  \end{gathered}
  \right\}\right.. 
\]
Since $f$ is regular, $I$ is a non-zero ideal of $R$. 

We claim that $\rho(I) \neq 0$. Indeed, 
suppose that $\rho(I) = 0$. Then the elimination theory 
(cf. \cite[Theorem~6]{KS}) tells us that there is a point 
$x = (x_1: \ldots : x_n) \in \PP^{N-1}(\widetilde{k})$
such that $\widetilde{F_i}(x,0)= 0$ and 
$\widetilde{G_j}(x,0)= 0$ for all $i$ and $j$ . 
Since $f$ satisfies the condition (ii), 
$\widetilde{F_i}(X,T)$ and $\widetilde{G_j}(X,T)$ are homogenization 
of $\widetilde{f_i}$ and $\widetilde{g_j}$, respectively. 
Then the existence of such an $x \in \PP^{N-1}(\widetilde{k})$ 
contradicts with the condition (iii). Hence we get the claim. 

Since $\rho(I) \neq 0$, there is an $r \in I$ such that $r \in 
R^{\times} = R \setminus M$. Then $I = R$, and we obtain 
Eqn.~\eqref{eqn:a:to:c}. 

(3) $\Longrightarrow$ (1): 
The assumption of (3) gives, as ideals in $\widetilde{k}[X]$, 
\[
 \left(X_1, \ldots, X_N\right)^\ell \subseteq 
 \left(\rho({F_1}(X,0)), \ldots, \rho({F_N}(X,0)), 
 \rho({G_1}(X,0)), \ldots, \rho({G_N}(X,0))\right). 
\]
Since $\rho({F_i}(X,0)) = \widetilde{F_i}(X,0)$ and 
$\rho({G_j}(X,0)) = \widetilde{G_j}(X,0)$, 
we obtain that $\widetilde{f}$ is regular. 

(2) $\Longrightarrow$ (3): We have only to put $T=0$. 
 
(3) $\Longrightarrow$ (2): It suffices to show that, for any 
$\alpha = 1, \ldots, N$, there are an integer $k \geq 1$ and 
polynomials $P_i(X,T)$, $Q_j(X,T)$ and $R(X,T)$ defined over $R$ such that 
\begin{equation}
\label{eqn:X:P:Q:R}
  X_\alpha^k = \sum_{i=1}^N P_i(X,T) F(X,T) + \sum_{j=1}^N Q_j(X,T) G_j(X,T) + 
  T R(X,T). 
\end{equation}
By the assumption of (iii), there are an integer $\ell \geq 1$ and 
polynomials $P_i(X)$, $Q_j(X)$ defined over $R$ such that 
\[
  X_\alpha^{\ell} 
  = \sum_{i=1}^N P_i(X) F(X,0) 
  + \sum_{j=1}^N Q_j(X) G_j(X,0). 
\]
We set $k := \ell$, $P_i(X,T) := P_i(X)$ and 
$Q_j(X,T) := Q_j(X)$. Then 
\[
  X_\alpha^k - \sum_{i=1}^N P_i(X,T) F(X,T) - \sum_{j=1}^N Q_j(X,T) G_j(X,T)
\]
is a polynomial in $R[X,T]$, which is divided by $T$. Hence there is a 
polynomial $R(X,T)$ in $R[X,T]$ satisfying Eqn.~\eqref{eqn:X:P:Q:R}. 
\QED

Suppose now that a regular polynomial automorphism $f$ 
has good reduction. By Proposition~\ref{prop:good:reduction:rephrasing}, 
for each $1 \leq i \leq N$, there are polynomials 
$P_{ij}(X)$ and $Q_{ij}(X)$ in $R[X]$ that satisfy \eqref{eqn:Pij:Qij}. 
Then the polynomial $R_i(X,T)$ in \eqref{eqn:pqr} is also defined 
over $R$. Then the constant $C$ in \eqref{eqn:C:def} is equal to $1$. 
This means that $\varepsilon = 1$ and $\delta = 1$ satisfy 
\eqref{eqn:e:d}. It follows that, when $f$ has good reduction,
$G_f$ and $\log^+\Vert\cdot\Vert$ are related simply. 

\begin{Proposition}
\label{prop:good:reduction:properties}
Suppose that $f$ has good reduction. Then 
\begin{enumerate}
\item[(1)]
$G_f(\cdot) \leq \log^+\Vert \cdot\Vert$ 
and $G_{f^{-1}}(\cdot) \leq \log^+\Vert \cdot\Vert$ 
on $\Aff^N(\Omega)$. 
\item[(2)]
$\log^+\Vert \cdot\Vert \leq G_f(\cdot)$ on ${V}_{1,1}^+$, 
$\log^+\Vert \cdot\Vert \leq G_{f^{-1}}(\cdot)$ on ${V}_{1,1}^-$, 
and $\Aff^N(\Omega) = {V}_{1,1}^+ \cup {V}_{1,1}^-$. 
\end{enumerate}
\end{Proposition}

\Proof
(1) Since the $f_i(X)$ are defined over $R$, 
in the proof of Lemma~\ref{lemma:G:above} 
we may take $r=1$, so that $c_f = 1$. Then 
$G_f(\cdot) \leq \log^+\Vert \cdot\Vert$ on $\Aff^N(\Omega)$. 
The estimate for $G_{f^{-1}}$  is similar. 

(2) 
Since $\varepsilon = 1$ and $\delta = 1$ satisfy 
\eqref{eqn:e:d}, Lemma~\ref{lemma:V:V-} gives $\Aff^N(\Omega) 
= {V}_{1,1}^+ \cup {V}_{1,1}^-$. 
The constant $c^+_{1,1}$ in Lemma~\ref{lemma:estimate:G} is 
equal to $0$, so we have 
$\log^+\Vert x\Vert \leq G_f(x)$ for all $x \in {V}_{1,1}^+$. 
The estimate for $G_{f^{-1}}$  is similar. 
\QED

\medskip
\section{Green functions for regular automorphisms over $\CC$}
\label{sec:Green:over:C}
In this section, we remark that 
the proof of Theorem~\ref{thm:local:II} 
gives a different proof (more explicit and without compactness arguments) of the corresponding 
estimates of Green functions over $\CC$ 

We write the usual absolute value of $\CC$ for $|\cdot|_\infty$, and we set 
$\Vert x \Vert_\infty := \max_i\{|x_i|_\infty\}$ for $x=(x_1 ,\ldots, x_N) 
\in \Aff^N(\CC)$. 

Let $f=(f_1, \ldots, f_N):\Aff^N\to\Aff^N$ be a regular polynomial automorphism 
of degree $d \geq 2$ defined over $\CC$.  Then the Green function for $f$ is defined 
by (see \cite[\S2]{Sib})
\begin{equation}
\label{eqn:Green:function:infty}
  \displaystyle{G_f(x) := \lim_{n\to+\infty}\frac{1}{d^n} 
  \log^+\Vert f^n(x)\Vert}
  \qquad 
  \text{for $x \in \Aff^N(\CC)$}. 
\end{equation}
Let $\Vert f\Vert_\infty$ be the maximum of the absolute values of all the coefficients of 
$f_i(X)$ for $1 \leq i \leq N$, and we set $c_{f, \infty} = \frac{1}{d-1}\log \max\left\{
\binom{N+d-1}{d}\Vert f\Vert_\infty, 1\right\}$. 
Note that $\binom{N+d-1}{d}$ is the number of monomials of 
degree $d$ in the ring of homogeneous polynomials in $N$ variables. 
Since 
\[
  \log^+\Vert f(x)\Vert \leq d \log^+\Vert f(x)\Vert 
  + \log \max\left\{\binom{N+d-1}{d}\Vert f\Vert_\infty, 1\right\}, 
\]
we get 
\begin{equation}
\label{eqn:Green:infinity}
 G_f(x) \leq \log^+\Vert x\Vert + c_{f, \infty}
 \qquad 
  \text{for any $x \in \Aff^N(\CC)$}. 
\end{equation}

Let $P_{ij}(X), Q_{ij}(X) \in \CC[X]$ and $R(X,T) \in \CC[X,T]$ 
be polynomials satisfying \eqref{eqn:pqr}. 
As before, we may and do assume that  the $P_{ij}(X)$, $Q_{ij}(X)$ and $R_i(X,T)$ are 
homogeneous polynomials with degree $m-d$, $m-d_-$ and $m-1$, 
respectively. We write $\Vert P\Vert_\infty$
for the maximum of the absolute values of all the coefficients of $P_{ij}(X)$ for $1 \leq i \leq N$ 
and $1 \leq j \leq N$, and we write $\Vert Q\Vert_\infty$, $\Vert R\Vert_\infty$ similarly. 
We set 
\[
  C_{\infty}^{\prime}
  = \max\left\{ \binom{N+ m-d -1}{m-d} \Vert P\Vert_\infty, 
  \binom{N+ m-d_- -1}{m-d_-} \Vert Q\Vert_\infty, 
  \binom{N+ m}{m-1} \Vert R\Vert_\infty, 1
  \right\}. 
\]
We put 
\[
  C_{\infty} = (2N+1) C_{\infty}^{\prime}. 
\]
We fix real numbers $\varepsilon>0, \delta>0$ satisfying \eqref{eqn:e:d} 
with $C_{\infty}$ in place of $C$. We define $N^{\pm}_{\delta, \varepsilon}$ 
and $V^{\pm}_{\delta, \varepsilon}$ by \eqref{eqn:n:v:+} and  \eqref{eqn:n:v:-}  
with $\CC$ in place of $\Omega$. 
Then, exactly as in Theorem~\ref{thm:local:II}, we have the following. 

\begin{Theorem}
\label{thm:local:II:infty}
$f: \Aff^N\to\Aff^N$ a regular polynomial automorphism 
over $\CC$. Then 
\begin{enumerate}
\item[(i)] $G_f(\cdot) \geq \log^+\Vert \cdot\Vert + c_{\delta, \varepsilon}^+$ 
on $V^{+}_{\delta, \varepsilon}$. 
\item[(ii)] $G_{f^{-1}}(\cdot) \geq \log^+\Vert \cdot\Vert + c_{\delta, \varepsilon}^-$ 
on $V^{-}_{\delta, \varepsilon}$. 
\item[(iii)] $V^{+}_{\delta, \varepsilon} \cup V^{-}_{\delta, \varepsilon} 
= \Aff^N(\Omega)$. 
\end{enumerate}
\end{Theorem}

\medskip
\section{Global theory of regular automorphisms}
\label{sec:global:theory}
\setcounter{equation}{0}
From this section, we turn our attention to regular automorphisms 
over a number field. 

Let $K$ be a number field, and $O_K$ its ring of integers. 
We fix an embedding 
$K \subset \overline{K}$ into an algebraic closure. 
Let $M_K$ be the set of absolute values on $K$. 
We extend the absolute values on $K$ to those on $\overline{K}$. 

Let $L$ be a finite extension 
field of $K$. For $x \in \Aff^N(L)$, we define 
\begin{equation}
\label{eqn:naive:height}
  h(x) = \sum_{v \in M_K} n_v \log^+\Vert x \Vert_v, 
\end{equation}
where $n_v = \frac{[L_v:K_v]}{[L:K]}$. 
This gives rise to the logarithmic Weil {\em height function}
\[
  h: \Aff^N(\overline{K}) \to \RR. 
\]
For more details on height functions, we refer 
the reader to \cite{BG, HS, La}. 

\smallskip
Let $f: \Aff^N \to \Aff^N$ be a regular polynomial automorphism 
over $\overline{K}$ (cf. Remark~\ref{rmk:regular:autom:def}). 
If the coefficients of $f$ are 
all defined over $K$, then we say that $f$ is a regular 
polynomial automorphism over $K$. 

\begin{Lemma}
\label{lemma:f-1:coeff}
If $f: \Aff^N \to \Aff^N$ is a regular polynomial automorphism 
over $K$, then the coefficients of $f^{-1}$ are 
also all defined over $K$. 
\end{Lemma}

\Proof
We take a finite Galois extension field $L$ of $K$ such that 
the coefficients of $f^{-1}$ are elements of $L$. For every 
$\sigma \in \Gal(L/K)$, the uniqueness of the inverse gives 
$(f^{-1})^{\sigma} = f^{-1}$. Thus the coefficients of $f^{-1}$ are 
in fact elements of $K$. 
\QED

\smallskip
In \cite{Ka}, we constructed (global) canonical height functions 
$h^+_f$ and $h^-_f$ for polynomial automorphisms $f$ over $K$, 
under the assumption that there exists a constant $c \geq 0$ 
such that 
\begin{equation}
\label{eqn:former:assumption:1}
  \frac{1}{d}h(f(x)) + \frac{1}{d_-}h(f^{-1}(x))
  \geq \left(1 + \frac{1}{d d_-}\right) h(x) -c  
\end{equation}
for all $x \in \Aff^N(\overline{K})$, 
where $d$ and $d_-$ denote the degrees of $f$ and $f^{-1}$. 
(We showed in {\em op.cit} that \eqref{eqn:former:assumption:1} 
holds for regular polynomial 
automorphisms in dimension $N = 2$ by a global method, i.e., 
a method using the effectiveness of a certain divisor 
on a certain rational surface.)

In the following, using properties of local Green functions 
studied in the previous sections, we will first construct in 
Theorem~\ref{thm:main} (global) canonical height functions $h^+_f$ and $h^-_f$ 
for regular polynomial automorphisms. 
Indeed, we will construct $h^+_f$ and $h^-_f$ 
as appropriate sums of local Green functions. 
Then, we show local versions of \eqref{eqn:former:assumption:1} for 
all places $v$, and summing them up   
we will obtain \eqref{eqn:former:assumption:1} for regular polynomial 
automorphisms in any dimension 
$N \geq 2$ in Theorem~\ref{thm:former:assumption}.

For a finite set $S$ of $M_K$ that contains all the Archimedean absolute values  of 
$K$, we denote by $(O_K)_S$ the ring of $S$-integers:
\[
  (O_K)_S = \{x \in K \mid \Vert x\Vert_v \leq 1 \quad\text{for all $v \not\in S$}\}. 
\]

\begin{Proposition}
\label{prop:good:reduction:are:ae}
Let $f: \Aff^N\to\Aff^N$ be a regular polynomial automorphism 
of degree $d \geq 2$ over a number field $K$. Then there 
exists a finite set $S$ of $M_K$ that contains all the Archimedean absolute values  of 
$K$ with the following property\textup{:} For all $v \not\in S$, $f$ induces the regular polynomial automorphism 
over $\overline{K_v}$ that has good reduction. 
\end{Proposition}

\Proof
We write $f= (f_1, \ldots, f_N)$ and let $F_i(X, T)  \in K[X, T]$ be the homogenization of 
$f_i$. Let $d_-$ denote the degree of $f^{-1} = (g_1, \ldots, g_N)$, and 
in virtue of Lemma~\ref{lemma:f-1:coeff} let 
$G_j(X, T) \in K[X, T]$ be the homogenization of $g_j$. Then there are 
an integer $m$ and homogeneous polynomials $P_{ij}(X) \in K[X]$ of degree $m-d$, 
 $Q_{ij}(X) \in K[X]$ of degree $m-d_-$, and  $R_i(X, T) \in K[X, T]$ of degree $m-1$ such 
 that \eqref{eqn:pqr} holds as polynomials in $K[X, T]$. 

We take a a finite set $S$ of $M_K$ that contains all the Archimedean absolute values  of 
$K$ with the following properties\textup{:} 
\begin{enumerate}
\item[(i)]
The coefficients of $F_i(X, T), G_j(X, T),  P_{ij}(X), Q_{ij}(X), R_i(X, T)$ are all in 
$(O_K)_S$. 
\item[(ii)]
For $v \not\in S$, we denote by $\rho_v: (O_K)_S \to \widetilde{k_v}$ the natural map, 
where $\widetilde{k_v}$ is the residue field of $(O_K)_v$. Then 
$\deg(f) = \deg(\rho_v(f))$ and $\deg(f^{-1}) = \deg(\rho_v(f^{-1}))$. 
\end{enumerate}
Then for any $v \not\in S$, 
$
  f \times_K \overline{K_v}: 
  \Aff^N_{\overline{K_v}} \to \Aff^N_{\overline{K_v}}
$
satisfies the properties (i) and (ii) of Definition~\ref{def:good:reduction} and (iii) of 
Proposition~\ref{prop:good:reduction:rephrasing}. Hence 
$ f \times_K \overline{K_v}$ has good reduction. 
\QED

\begin{Theorem}
\label{thm:main}
Let $f: \Aff^N\to\Aff^N$ be a regular polynomial automorphism 
of degree $d \geq 2$ over a number field $K$. Let $d_- \geq 2$ 
denote the degree of $f^{-1}$. 
\begin{enumerate}
\item[(1)] 
For all $x \in \Aff^N(\overline{K})$, 
the limits 
\[
  \lim_{n\to+\infty}\frac{1}{d^n} h(f^n(x)) 
  \quad\text{and}\quad
  \lim_{n\to+\infty}\frac{1}{d_-^n} h(f^{-n}(x))
\]
exist. We write $\widehat{h}^+_f(x)$ and $\widehat{h}^-_{f}(x)$ 
for the limits, respectively.
\item[(2)]\textup{(}Global-to-local decomposition\textup{)}
For each place $v \in M_K$, 
let $G_{f,v}$ and $G_{f^{-1},v}$ be the Green functions 
for $f$ and $f^{-1}$ at $v$, respectively. 
Then, for all $x \in \Aff^N(\overline{K})$, one has 
\[
  \widehat{h}^+_f(x) = \sum_{v \in M_K} n_v G_{f,v}(x) 
  \quad\text{and}\quad
  \widehat{h}^-_{f}(x) = \sum_{v \in M_K} n_v G_{f^{-1},v}(x) 
\]
\item[(3)]
We define $\widehat{h}_f: \Aff^N(\overline{K})\to\RR$ by 
\[
  \widehat{h}_f := \widehat{h}^+_f + \widehat{h}^-_{f}. 
\]
Then $\widehat{h}_f$ satisfies the following two conditions. 
\begin{enumerate}
\item[(i)]
$\displaystyle 
\frac{1}{d} \widehat{h}_f\circ f + 
\frac{1}{d_-} \widehat{h}_f\circ f^{-1} 
= \left(1+\frac{1}{d d_-}\right) \widehat{h}_f$ 
on $\Aff^N(\overline{K})$. 
\item[(ii)]
$\widehat{h}_f \gg\ll h$ on $\Aff^N(\overline{K})$. 
\end{enumerate}
\item[(4)]
The function $\widehat{h}_f$ has the following uniqueness property\textup{:} 
If $h^{\prime}: \Aff^N(\overline{K})\to\RR$ is a function 
satisfying the condition \textup{(3-i)}  
such that $h^{\prime} = \widehat{h}_f +O(1)$, 
then $h^{\prime} = \widehat{h}_f$. 
\item[(5)]
The functions $\widehat{h}^+_f$, $\widehat{h}^-_f$ 
and $\widehat{h}_f$ are non-negative. Further, 
for $x \in \Aff^N(\overline{K})$, we have 
\[
  \widehat{h}_f(x) = 0 
  \;\Longleftrightarrow\;
  \widehat{h}^+_f(x) = 0 
  \;\Longleftrightarrow\;
  \widehat{h}^-_f(x) = 0 
  \;\Longleftrightarrow\;
  \text{$x$ is $f$-periodic.}
\]
\end{enumerate}
\end{Theorem}

\Proof
For each $v \in M_K$, we have estimates of Green functions 
for $f$ at $v$ as in Lemma~\ref{lemma:G:above} and 
Lemma~\ref{lemma:estimate:G}. We use the suffix $v$ when we work over the absolute 
value $v \in M_K$. For example, Green function for $f$ at $v$ is denoted by $G_{f, v}$ 
and constants $c_f, c^\pm_{\varepsilon, \delta}$ in Lemma~\ref{lemma:G:above}, 
Lemma~\ref{lemma:estimate:G} and \eqref{eqn:estimate:G-} are denoted by 
$c_{f, v}, c^\pm_{\varepsilon, \delta, v}$ respectively. 

Let $S$ be a finite subset of $M_K$ in Proposition~\ref{prop:good:reduction:are:ae}. 

(1)(2) We fix $x \in \Aff^N(\overline{K})$. 
We will show the existence of 
$h_f^+(x)$ and the decomposition 
$h_f^+(x) = \sum_{v \in M_K} n_v G_{f,v}(x)$. 
The existence and decomposition for $h_f^-(x)$ are shown similarly. 

For $v \in M_K$ and $n \geq 0$, we set 
\[
  G^+_{v, n}(x) := \frac{1}{d^n} \log^+ \Vert f^n(x)\Vert_v. 
\]
Then we have 
\begin{itemize}
\item 
$0 \leq G^+_{v, n}(x) \leq \log^+\Vert x \Vert_v + c_{f, v}$ 
for all $v \in M_K$ and $n \geq 0$ from Lemma~\ref{lemma:G:above} (or its proof)
and  Eqn.~\eqref{eqn:Green:infinity}. 
\item 
$\lim_{n \to +\infty}  G^+_{v, n}(x) = G_{f, v}(x)$ 
from Definition~\ref{def:Green:function} and Eqn.~\eqref{eqn:Green:function:infty}. 
\item
$\frac{1}{d^n} h(f^n(x))  = \sum_{v \in M_K} n_v G^+_{v, n}(x)$ 
from Eqn.~\eqref{eqn:naive:height}.  
\item 
We may take $c_{f,v} = 0$ for any $v \not\in S$ from 
Proposition~\ref{prop:good:reduction:properties} and 
Proposition~\ref{prop:good:reduction:are:ae}.   
\item
$\sum_{v \in M_K} n_v (\log^+\Vert x \Vert_v + c_{f,v}) 
= h(x) + \sum_{v \in S} n_v c_{f,v} < +\infty$. 
\end{itemize}
Lebesgue's dominant convergence theorem then implies that 
$\sum_{v \in M_K} n_v G^+_{v, n}(x)$ 
converges as $n\to+\infty$ and that   
\begin{align*}
  \lim_{n\to+\infty} \frac{1}{d^n} h(f^n(x))  
  & = \lim_{n\to+\infty} \sum_{v \in M_K} n_v G^+_{v, n}(x) \\
  & = \sum_{v \in M_K} \lim_{n\to+\infty} n_v G^+_{v, n}(x) 
    = \sum_{v \in M_K} n_v G_{f, v}(x). 
\end{align*}
This completes the proof of (1)(2). 

\smallskip
(3)(4)(5) 
First we have 
\begin{align}
\label{eqn:3:1}
  \widehat{h}_f(x) 
  & =  \sum_{v \in M_K} n_v G_{f, v}(x) + 
  \sum_{v \in M_K} n_v G_{f^{-1}, v}(x) \\
\notag    
  & \leq 
  \sum_{v \in M_K} n_v (2 \log^+ \Vert x\Vert_v 
  + c_{f,v} + c_{f^{-1},v}) 
  = 2 \widehat{h}_{nv}(x) + 
  \sum_{v \in S} n_v (c_{f,v} + c_{f^{-1},v}). 
\end{align}

On the other hand, we have 
\begin{itemize}
\item 
$\min\{c^+_{\varepsilon, \delta, v}, 
c^-_{\varepsilon, \delta, v}\} 
+ \log^+ \Vert x\Vert\leq G_{f, v}(x) + G_{f^{-1}, v}(x)$ 
from Lemma~\ref{lemma:estimate:G}, 
Eqn.~\eqref{eqn:estimate:G-} and Theorem~\ref{thm:local:II:infty}. 
\item 
For any $v \not\in S$, we may take $\varepsilon=1, \delta=1$ and 
$\min\{c^+_{1, 1, v}, 
c^-_{1, 1, v}\} = 0$ from 
Proposition~\ref{prop:good:reduction:properties} and 
Proposition~\ref{prop:good:reduction:are:ae}.   
\end{itemize}
Then 
\begin{align}
\label{eqn:3:2}
  \widehat{h}_f(x) 
  & =  \sum_{v \in M_K} n_v G_{f ,v}(x) + 
  \sum_{v \in M_K} n_v G_{f^{-1}, v}(x) \\
\notag  
  & \geq  \sum_{v \in M_K} n_v \left( \log^+ \Vert x\Vert_v
  + \min\{c^+_{\varepsilon, \delta, v}, 
  c^-_{\varepsilon, \delta, v}\}\right) 
  = 
  \widehat{h}_{nv}(x) + \sum_{v \in S} n_v \min\{c^+_{\varepsilon, \delta, v}, 
  c^-_{\varepsilon, \delta, v}\}. 
\end{align}
Eqns.~\eqref{eqn:3:1} and \eqref{eqn:3:2} give (3)(ii). 
For the rest of the proof, 
see \cite[Theorem~4.2(2)(3)(4)]{Ka}.
\QED

\begin{Remark}
Theorem~\ref{thm:main}(1) shows that 
$\{   
\frac{1}{d^n} h(f^n(x))\}_{n=0}^{+\infty}$ and $\{   
\frac{1}{d_-^n} h(f^{-n}(x))\}_{n=0}^{+\infty}$
are convergent sequences, which gives an improvement of \cite{Ka}, 
since we replace $\limsup$ by $\lim$ in the definition 
of $\widehat{h}^\pm_f$.    
\end{Remark}

We now introduce another function 
\begin{equation}
\label{eqn:another:can:ht:fcn}
  \widetilde{h}_f(x) := 
  \sum_{v \in M_K} n_v 
  \max\{G_{f,v}(x) , G_{f^{-1},v}(x) \}
\end{equation}
for $x \in \Aff^N(\overline{K})$. 
The next proposition shows that  $\widetilde{h}_f$ 
also behaves well relative to $f$.  

\begin{Proposition}
\label{prop:tilde:G}
\begin{enumerate}
\item
$\widetilde{h}_f = h + O(1)$ on 
$\Aff^N(\overline{K})$. 
\item
For $x \in \Aff^N(\overline{K})$, we have 
$\widetilde{h}_f(x) = 0$ if and only if $\widehat{h}_f(x) = 0$.
\end{enumerate}
\end{Proposition}

\Proof
(1) We use the notation of the proof of Theorem~\ref{thm:main}. 
By Lemma~\ref{lemma:G:above}, Eqn.~\eqref{eqn:Green:infinity},  
Lemma~\ref{lemma:estimate:G}, Eqn.~\eqref{eqn:estimate:G-} and 
Theorem~\ref{thm:local:II:infty}, 
we have 
\[
 \log^+ \Vert x\Vert_v
  + \min\{c^+_{\varepsilon, \delta, v}, c^-_{\varepsilon, \delta, v}\}
  \leq 
  \max\{G_{f,v}(x) , G_{f^{-1},v}(x)\}  
  \leq 
\log^+ \Vert x\Vert_v + \max\{
  c_{f,v}, c_{f^{-1}, v}\}.  
\]
Summing up over all places $v$, we get 
\[
  h(x) + \sum_{v \in M_K} n_v \min\{c^+_{\varepsilon, \delta, v}, c^-_{\varepsilon, \delta, v}\}
  \leq 
  \widetilde{h}_f (x)
  \leq 
  h(x) + \sum_{v \in M_K} n_v \max\{
  c_{f,v}, c_{f^{-1}, v}\}.
\]
Since   we have 
$c_{f,v} = c_{f^{-1}, v} = c^+_{\varepsilon, \delta, v} = c^-_{\varepsilon, \delta, v} =0$ 
except for finitely many $v$ (indeed for every $v \not\in S$), 
this gives the assertion. 

(2)
Since $G_{f,v}$ and $G_{f^{-1},v}$ are non-negative functions, 
we see that $\widetilde{h}_f(x) = 0$  if and only if 
$G_{f,v}(x) = G_{f^{-1},v}(x) = 0$ if and only if $\widehat{h}_f(x) = 0$. 
\QED

\medskip
\section{Arithmetic properties of regular polynomial automorphisms}
\label{sec:applications}
\setcounter{equation}{0}
In this section, we give some applications of local and 
global canonical height functions. 
The first application is the following theorem on 
the usual height function (see \cite[\S4]{Ka}, \cite[Conjecture~3]{Sil}, \cite[Conjecture~7.18]{SilBook}), 
which is independently obtained by Lee \cite{Lee} via a different method (via a global method 
based on the effectiveness of a certain divisor as in the case of $N=2$ in \cite{Ka}). 

\begin{Theorem}
\label{thm:former:assumption}
Let $f: \Aff^N\to\Aff^N$ be a regular polynomial automorphism 
over a number field $K$. Let $d$ and $d_-$ be the degrees 
of $f$ and $f^{-1}$. Then we have the following. 
\begin{enumerate}
\item[(1)]
There exists a constant $c \geq 0$ 
such that 
\[
  \frac{1}{d}h(f(x)) + \frac{1}{d_-}h(f^{-1}(x))
  \geq \left(1 + \frac{1}{d d_-}\right) h(x) -c  
\]
for all $x \in \Aff^N(\overline{K})$.  
\item[(2)]
One has  
\begin{equation}
\label{eqn:former:assumption:2}
  \liminf_{\substack{
  x \in \Aff^N(\overline{K})\\ 
  h(x) \to \infty}} 
  \frac{
  \frac{1}{d}h(f(x)) + \frac{1}{d_-}h(f^{-1}(x))}{h(x)} 
  = 1 + \frac{1}{d d_-}. 
\end{equation}
\end{enumerate}
\end{Theorem}

\Proof
(1) 
We set 
\[
  \widetilde{G}_{f, v} := \max\{G_{f,v}, G_{f^{-1},v}\}. 
\]

\begin{Claim}
\label{claim:tilde:G}
For all $x \in \Aff^N(\overline{K})$, we have 
\begin{equation}
\label{eqn:former:assumption:G}
  \frac{1}{d}\widetilde{G}_{f,v}(f(x)) + \frac{1}{d_-}\widetilde{G}_{f,v}(f^{-1}(x))
  \geq \left(1 + \frac{1}{d d_-}\right) \widetilde{G}_{f,v}(x). 
\end{equation}
\end{Claim}

We first show that Claim~\ref{claim:tilde:G} implies (1). Indeed, we assume 
Claim~\ref{claim:tilde:G}. Then summing up over all $v$, we have 
\begin{equation}
\label{eqn:former:assumption}
  \frac{1}{d}\widetilde{h}(f(x)) + \frac{1}{d_-}\widetilde{h}(f^{-1}(x))
  \geq \left(1 + \frac{1}{d d_-}\right) \widetilde{h}(x).   
\end{equation}
Since $\widetilde{h}_f = h +O(1)$ by Proposition~\ref{prop:tilde:G}(1), 
Eqn.~\eqref{eqn:former:assumption} yields (1). 

Recall that $G_{f,v}$ and $G_{f^{-1},v}$ are non-negative. 
We will show \eqref{eqn:former:assumption:G} according to 
the  four cases:  
$\frac{1}{d d_-} G_{f,v}(x) \gtrless G_{f^{-1},v}(x)$ 
and $\frac{1}{d d_-} G_{f^{-1},v}(x) \gtrless G_{f,v}(x)$. 
The first case treats $(\geq, \leq)$, and then $(\leq, \leq)$, $(\geq, \geq)$ and 
$(\leq, \geq)$ in this order. 

{\bf Case 1:}\quad Suppose that 
$\frac{1}{d d_-} G_{f,v}(x) \geq G_{f^{-1},v}(x)$. 
It follows that $G_{f,v}(x)\geq d d_- G_{f^{-1},v}(x) \geq \frac{1}{d d_-}G_{f^{-1},v}(x)$. 
In this case, we have 
\begin{align*}
  \widetilde{G}_{f,v}(f(x)) 
  & = \max\{G_{f,v}(f(x)), G_{f^{-1},v}(f(x))\}  \\
  & = \max\left\{d G_{f,v}(x), \frac{1}{d_-}G_{f^{-1},v}(x)\right\} 
  = d G_{f,v}(x),  
\end{align*}
and similarly 
\begin{align*}
  \widetilde{G}_{f,v}(f^{-1}(x)) 
  & = \max\{G_{f,v}(f^{-1}(x)), G_{f^{-1},v}(f^{-1}(x))\} \\
  & = \max\left\{\frac{1}{d} G_{f,v}(x), d_- G_{f^{-1},v}(x)\right\} 
  = \frac{1}{d} G_{f,v}(x). 
\end{align*}
Further, the assumption $\frac{1}{d d_-} G_{f,v}(x) 
\geq G_{f^{-1},v}(x)$ implies that $G_{f,v}(x) \geq G_{f^{-1},v}(x)$.  
Thus $\widetilde{G}_{f, v}(x) =G_{f,v}(x)$. Putting these together, we get  
\[
  \text{Left-hand side of \eqref{eqn:former:assumption:G}} 
  = \left(1 + \frac{1}{d d_-}\right) G_{f,v}(x)
  = \text{Right-hand 
  side of \eqref{eqn:former:assumption:G}}.
\]

{\bf Case 2:}\quad Suppose that 
$G_{f,v}(x) \geq \frac{1}{d d_-} G_{f^{-1},v}(x)$ 
and $\frac{1}{d d_-} G_{f,v}(x) \leq G_{f^{-1},v}(x)$. 
In this case, we have 
$ \widetilde{G}_{f,v}(f(x)) = d G_{f,v}(x)$ 
and 
$ \widetilde{G}_{f,v}(f^{-1}(x))  = d_- G_{f^{-1},v}(x)$. 

{\bf Subcase 2-1:}\quad Suppose that 
$G_{f,v}(x)  \geq G_{f^{-1},v}(x)$. Then, 
using $\frac{1}{d d_-} G_{f,v}(x)  \leq G_{f^{-1},v}(x)$, we have 
\[
  \text{Left-hand side of \eqref{eqn:former:assumption:G}} 
  = G_{f,v}(x) +  G_{f^{-1},v}(x) \\
  \geq \left(1 + \frac{1}{d d_-}\right) G_{f,v}(x)
  = \text{Right-hand 
  side of \eqref{eqn:former:assumption:G}}.
\]

{\bf Subcase 2-2:}\quad Suppose that 
$G_{f^{-1},v}(x) \geq G_{f,v}(x)$. Then, 
using  $G_{f,v}(x) \geq \frac{1}{d d_-} G_{f^{-1},v}(x)$, 
\[
  \text{Left-hand side of \eqref{eqn:former:assumption:G}} 
  = G_{f,v}(x) +  G_{f^{-1},v}(x)  \\
  \geq \left(1 + \frac{1}{d d_-}\right) G_{f^{-1},v}(x)
  = \text{Right-hand 
  side of \eqref{eqn:former:assumption:G}}.
\]

{\bf Case 3:}\quad Suppose that 
$G_{f,v}(x) \leq \frac{1}{d d_-} G_{f^{-1},v}(x)$ 
and $\frac{1}{d d_-} G_{f,v}(x) \geq G_{f^{-1},v}(x)$.
Then $G_{f^{-1},v}(x) \leq \frac{1}{d d_-} G_{f,v}(x)
\leq \left(\frac{1}{d d_-} \right)^2 G_{f^{-1},v}(x)$. 
It follows that $G_{f,v}(x) = G_{f^{-1},v}(x) = 0$. 
Then 
\[
  \widetilde{G}_{f,v}(f(x)) = \max\left\{
  G_{f,v}(f(x)),  G_{f^{-1},v}(f(x))
  \right\} 
  = 
   \max\left\{
  d G_{f,v}(x),  \frac{1}{d_-}G_{f^{-1},v}(x)
  \right\} 
  =0, 
\]
and similarly we have $\widetilde{G}_{f,v}(f^{-1}(x)) = 0$ and 
$\widetilde{G}_{f,v}(x) =0$. 
We get 
\[
  \text{Left-hand side of \eqref{eqn:former:assumption:G}} 
  = 0
  = \text{Right-hand 
  side of \eqref{eqn:former:assumption:G}}.
\]

{\bf Case 4:}\quad Suppose that 
$G_{f,v}(x)  \leq \frac{1}{d d_-} G_{f^{-1},v}(x)$. It follows that 
$G_{f^{-1},v}(x) \geq d d_- G_{f,v}(x)  \geq \frac{1}{d d_-}G_{f,v}(x) $. 
We can show \eqref{eqn:former:assumption:G} as in Case~1, 
exchanging the roles of $G_{f,v}$ and $G_{f^{-1},v}$. 

This completes the proof of Claim~\ref{claim:tilde:G}, hence 
the proof of Theorem~\ref{thm:former:assumption}(1). 

\medskip
(2) From (1), we obtain 
\[
  \liminf_{\substack{
  x \in \Aff^N(\overline{K})\\ 
  h(x) \to \infty}} 
  \frac{
  \frac{1}{d}h(f(x)) + \frac{1}{d_-}h(f^{-1}(x))}{h(x)} 
  \geq 1 + \frac{1}{d d_-}. 
\]
On the other hand, it is shown in \cite[Proposition~4.4]{Ka} that, 
for any polynomial automorphism $f:\Aff^N\to\Aff^N$, 
one has 
\begin{equation}
  \liminf_{\substack{
  x \in \Aff^N(\overline{K})\\ 
  h(x) \to \infty}} 
  \frac{
  \frac{1}{d}h(f(x)) + \frac{1}{d_-}h(f^{-1}(x))}{h(x)} 
  \leq 1 + \frac{1}{d d_-}.   
\end{equation}
Combining these two inequalities gives the assertion. 
\QED

\begin{Remark}
It is shown in \cite[Theorem~4.4]{Ka} that the equality  
\eqref{eqn:former:assumption:2} holds in dimension $N=2$ 
for regular polynomial automorphisms. 
Theorem~\ref{thm:former:assumption}(2) asserts that 
the equality holds  in any dimension $N \geq 2$ for 
regular polynomial automorphisms. 
\end{Remark}

Theorem~\ref{thm:main} recovers the following 
theorem of Marcello \cite{Ma1} on $f$-periodic points. 

\begin{Corollary}[\cite{Ma1}]
\label{cor:Mar}
Let $f: \Aff^N\to\Aff^N$ be a regular polynomial automorphism 
over a number field $K$. Then 
the set of $f$-periodic points in $\Aff^N(\overline{K})$ 
is a set of bounded height. In particular, for any integer 
$D \geq 1$, the set 
\[
  \left\{
  x \in \Aff^N(\overline{K}) \;\left\vert\;
  \text{$x$ is $f$-periodic},\; [K(x): K] \leq D \right. 
  \right\}
\]
is finite. 
\end{Corollary}

\Proof
By Theorem~\ref{thm:main}(3-ii) and (5), 
$\widehat{h}_{f}$ satisfies $\widehat{h}_{f} \gg\ll h$, 
and a point $x \in \Aff^N(\overline{K})$ is $f$-periodic 
if and only if $\widehat{h}_{f}(x) = 0$. Thus we get the assertion. 
\QED

For a non $f$-periodic point $x$, let 
$O_f(x) := \{f^n(x) \mid n \in \ZZ\}$ denote the $f$-orbit of $x$. 
We define the canonical height of the orbit $O_f(x)$ by 
\begin{equation}
\label{eqn:can:ht:orbit}
  \widehat{h}_f(O_f(x)) 
  = \frac{\log \widehat{h}^+(x)}{\log d} 
  + \frac{\log \widehat{h}^-(x)}{\log d_-} , 
\end{equation}
whose value depends only on the orbit $O_f(x)$ 
and not the particular choice of the point 
$x$ in the orbit by Theorem~\ref{thm:main}.  
The next corollary gives a refinement of 
 \cite[Corollaire~B]{Ma2}.

\begin{Corollary}
\label{cor:refined}
Let $f: \Aff^N\to\Aff^N$ be a regular polynomial automorphism 
over a number field $K$. 
Let $d$ and $d_-$ be the degrees of $f$ and $f^{-1}$. Then 
for any infinite orbit $O_f(x)$, 
\[
  \#\{y \in O_f(x) \mid h(y) \leq T\}
  = 
  \left(\frac{1}{\log d} + \frac{1}{\log d_-} \right)\log T 
  - \widehat{h}(O_f(x)) + O(1)
\]
as $T \to +\infty$. Here the $O(1)$ bound depends only $f$, 
independent of the orbit $O_f(x)$. 
\end{Corollary}

\Proof
Since $f$ satisfies \eqref{eqn:former:assumption}, 
we apply \cite[Theorem~5.2]{Ka}. 
\QED

In the rest of this section, we consider some global-to-local 
arithmetic properties. 
Suppose that $f$ is a regular polynomial automorphism. 
By Theorem~\ref{thm:main}(2)(5), 
$x \in \Aff^N(\overline{K})$ is $f$-periodic 
if and only if $G_{f,v}(x) = 0$ for all $v \in M_K$. 
By Theorem~\ref{thm:G:K} for non-Archimedean $v$ and \cite[\S2]{Sib} for Archimedean $v$, 
$G_{f,v}(x)=0$ is equivalent to  
$\{f^n(x)\}_{n=0}^{+\infty}$ being bounded 
with respect to $\Vert\cdot\Vert_v$. 
Thus we see that $x \in \Aff^N(\overline{K})$ is $f$-periodic 
if and only if $\{f^n(x)\}_{n=0}^{+\infty}$ is bounded 
with respect to $\Vert\cdot\Vert_v$ for all $v \in M_K$. 

This actually holds for any polynomial map $f$ (cf. 
\cite[Corollary~6.3]{CG} for $N=1$). 

\begin{Proposition}
Let $f: \Aff^N\to\Aff^N$ be a polynomial map over a number field $K$. 
For $x \in \Aff^N(\overline{K})$, the following is equivalent.  
\begin{enumerate}
\item[(i)]
$x$ is $f$-periodic. 
\item[(ii)]
For every $v \in M_K$, 
$\{f^n(x)\}_{n=0}^{+\infty}$ is bounded 
with respect to the $v$-adic topology. 
\end{enumerate}
\end{Proposition}

\Proof
Taking a finite extension field of $K$ over which $x$ is defined if necessary, 
we may assume that $x$ is defined over $K$. 
It is obvious that (i) implies (ii). We assume (ii) and 
show (i). We take a finite subset $S$ of $M_K$ containing the set of all Archimedean absolute 
values such that 
$x$ and $f$ is defined over $(O_K)_S$. Then for any $v \not\in S$, 
we have 
\[
  \Vert f^n(x)\Vert_v \leq 1 
  \qquad\text{for all $n \geq 0$}.
\]
Since we assume (ii), 
there is a constant $C_v$ for each $v \in S$ such that 
\[
  \Vert f^n(x)\Vert_v \leq C_v 
  \qquad\text{for all $n \geq 0$}. 
\] 
Then we have 
\[
  h(f^n(x)) 
  = \sum_{v \in M_K} n_v \log^+\Vert f^n(x)\Vert 
  \leq \sum_{v \in S} n_v C_v
  \qquad\text{for all $n \geq 0$}. 
\]
Then 
\[
  \{ f^n(x) \mid n \geq 0\} 
  \quad\subseteq\quad  
  \{ y \in \Aff^N(K) \mid h(y) \leq\sum_{v \in S} n_v C_v \}.  
\]
Since the latter set is finite, the set $\{f^n(x)\}_{n \geq 0}$ is finite. 
Hence $x$ is $f$-periodic.  
\QED

\end{document}